\documentclass[10pt,A4paper]{article}

\usepackage{tikz} 
\usepackage{amsmath}
\usepackage{mathrsfs} 
\usepackage{amsfonts} 
\usepackage{pgfplots}

 \usepackage{amssymb}

 \usepackage[colorlinks, linkcolor=blue, anchorcolor=red, citecolor=blue]{hyperref}

 \usepackage{appendix}

\def\ar{\!\!\!&} 

\def\proof{\noindent{\it Proof.~~}} 
\def\qed{\hfill$\Box$\medskip} 
 \usepackage[top=2cm,bottom=2cm, outer=2.3cm, inner=2.3cm]{geometry}

\newtheorem{theorem}{Theorem}[section]

\newtheorem{lemma}[theorem]{Lemma}
\newtheorem{corollary}[theorem]{Corollary}
\newtheorem{proposition}[theorem]{Proposition}
\newtheorem{remark}[theorem]{Remark}
\newtheorem{condition}[theorem]{Condition}
\newtheorem{example}[theorem]{Example}

\def\beqlb{\begin{eqnarray}}\def\eeqlb{\end{eqnarray}} 
\def\beqnn{\begin{eqnarray*}}\def\eeqnn{\end{eqnarray*}} 

\def\ar{\!\!\!&}

\def\proof{\noindent{\it Proof.~~}}\def\qed{\hfill$\Box$\medskip}

\begin{document}
 \title{\bf\Large Functional Limit Theorems for Marked Hawkes Point Measures\thanks{Financial support by the Alexander-von-Humboldt Foundation is gratefully acknowledged.}}
 
 \author{Ulrich Horst\footnote{Department of Mathematics and School of Business and Economics,  Humboldt-Universit\"at zu Berlin, Unter den Linden 6, 10099 Berlin; email: horst@math.hu-berlin.de}\quad\  and\quad Wei Xu\footnote{Department of Mathematics, Humboldt-Universit\"at zu Berlin, Unter den Linden 6, 10099 Berlin; email: xuwei@math.hu-berlin.de}}
 \maketitle


  \begin{abstract}

This paper establishes a functional law of large numbers and a functional central limit theorem for marked Hawkes point measures and their corresponding shot noise processes. We prove that the normalized random measure can be approximated in distribution by the sum of a Gaussian wihte noise process plus an appropriate lifting map of a correlated one-dimensional Brownian motion. The Brownian results from the self-exiting arrivals of events. We apply our limit theorems for Hawkes point measures to analyze the population dynamics of budding microbes in a host.    

  \end{abstract}



  \section{Introduction}
\setcounter{equation}{0}
\medskip

Let $(\Omega,\mathscr{F},\mathbf{P})$ be a complete probability space endowed with filtration $\{\mathscr{F}_t:t\geq 0\}$ that satisfies the usual hypotheses and $\mathbb{U}$ be a Lusin topological space endowed with the Borel $\sigma$-algebra $\mathscr{U}$. Let $\{\tau_k:k=1,2\cdots\}$  be a sequence of increasing, $(\mathscr{F}_t)$-adaptable random times and $\{ \xi_k:k=1,2,\cdots \}$ be  a sequence of i.i.d. $\mathbb{U}$-valued random variables with distribution $\nu_H(du)$. 
We assume that $\xi_k$ is independent of  $\{\tau_j:j=1,\cdots,k \}$ for any $k\geq 0$.
In terms of these sequences we define the $(\mathscr{F}_t)$-random point measure 
\beqlb\label{eqn2.03}
N_H(ds,du):= \sum_{k=1}^\infty \mathbf{1}_{\{\tau_k\in ds, \xi_k\in du  \}}
\eeqlb
on $(0,\infty)\times \mathbb{U}$. We say $N_H(ds,du)$  is a \textit{marked Hawkes point measure} if the embedded point process $\{N_t:t\geq 0  \}$ defined by 
\beqlb
	N_t:= N_H((0,t],\mathbb{U})
\eeqlb
admits an $(\mathscr{F}_t)$-intensity $\{ Z(t):t\geq 0 \}$ of the form 
\beqlb\label{eqn1.02}
Z(t)\ar=\ar \mu(t) + \sum_{k=1}^{N_t} \phi(t-\tau_k, \xi_k), \quad t\geq 0,
\eeqlb
for some nonnegative, locally integrable, $(\mathscr{F}_t)$-progressive  \textit{exogenous intensity} $\{\mu(t):t\geq 0\}$,  and some \textsl{kernel} $\phi: \mathbb{R}_+\times \mathbb{U} \to [0,\infty)$. 
%
We call $N_H(ds,du)$ a marked Hawkes point measures with homogeneous immigration if the exogenous intensity $\{\mu(t):t\geq 0\}$ is driven by an independent marked Poisson point measure 
 \beqlb\label{eqn2.01}
 N_I(ds,du)\ar:=\ar \sum_{k=1}^\infty \mathbf{1}_{\{\sigma_k\in ds, \eta_k\in du  \}}
 \eeqlb
on $(0,\infty)\times \mathbb{U}$ where the immigration times  $\{\sigma_j:j=1,2,\cdots \}$ are described by a Poisson point process 
\beqlb
	N'_t:=N_I([0,t),\mathbb{U})
\eeqlb
with some rate $\lambda_I$, and the marks are described by an i.i.d. sequence $\{ \eta_k:k=1,2,\cdots \}$ of $\mathbb{U}$-valued random variables with distribution $\nu_I(du)$ that is independent of $N'$. Specifically, we say that $N_H(ds,du)$ is a marked Hawkes point measures with homogeneous immigration if $\{\mu(t):t\geq 0\}$  admits the representation 
\beqlb\label{eqn1.03}
\mu(t)\ar:=\ar \mu_0(t) + \sum_{k=1}^{N'_t} \phi(t-\sigma_k, \eta_k),\quad t\geq 0,
\eeqlb
for some $\mathscr{F}_0$-measurable, nonnegative functional-valued random variable  $\{  \mu_0(t) :t\geq 0 \}$ that describes the impact of events prior to time $0$ on the arrival of future events. 


	
Marked Hawkes point measures with homogeneous immigration contain several important processes as special cases. If the exogenous intensity $\{\mu(t):t\geq 0\}$ equals some deterministic constant and the mark space is finite, then the measure $N_H(ds,du)$ reduces to a multivariate Hawkes process with common intensity; general multivariate Hawkes processes correspond to multi-dimensional Hawkes point measures whose mark spaces contain a single element. For constant exogenous intensities the embedded point process $\{N_t:t\geq 0  \}$ reduces to a marked Hawkes  process; see  \cite{BremaudMassoulie2002, Bremaud2002}. If, in addition, the kernel $\phi$ is independent of the mark $\xi$, then it reduces to a standard Hawkes process. 
First introduced in \cite{Hawkes1971a, Hawkes1971b} Hawkes processes have long become a powerful tool to model a variety of phenomena in science and finance; we refer to  \cite{BacryMastromatteoMuzy2015, Bordenave2007} for reviews on Hawkes processes and their applications. Compared to  Hawkes processes, marked Hawkes processes are ``individual based". The marks can be considered as the characteristics of an event, with events with different characteristics having different impacts on the arrival of future events. For instant, earthquakes of difference magnitudes have different effects on the arrivals of the future earthquakes \cite{Ogata1986}; in limit order markets, market orders of different sizes have different impacts on arrivals of future orders \cite{BacryIugaLasnierLehalle2015, Chavez-Demoulin2012}; in electricity markets price spikes of different sizes have different impacts on the occurrence of future spikes \cite{ClementsHerreraHurn2015}.   

In this paper, we establish functional limit theorems for marked Hawkes point measures with homogeneous immigration. As a byproduct we obtain a novel functional CLT for marked Hawkes processes. 
A {functional} CLT for multivariate Hawkes processes has been established by Bacry et al. \cite{BacryDelattreHoffmannMuzy2013}. A CLT for marked Hawkes processes has been proved by Karabash and Zhu \cite{KarabashZhu2015}; it can also be derived from the large deviation principle proved in Stabile and Torriai \cite{StabileTorrisi2010} and Zhu  \cite{Zhu2013}. In a recent paper Gao and Zhu \cite{GaoZhu2018, GaoZhu2018b} established a functional CLT and a large deviation principle for Hawkes processes with exponential kernel and large initial intensity. For the nearly unstable Hawkes process, Jaisson and Rosenbaum \cite{JaissonRosenbaum2015} proved that the rescaled intensity  converges weakly to a Feller diffusion and that the rescaled point process converges weakly to the integrated diffusion.

In order to prove our {functional} limit theorems we first establish a functional CLT for the embedded marked Hawkes process $\{N_t:t\geq 0  \}$. 
The key is to analyze the covariance function $\{ {\rm Cov}(Z(t), Z(t+s)):s,t\geq 0 \}$ of the intensity process $\{Z(t):t\geq 0\}$.
Compared to Hawkes processes, the analysis of the intensity of marked Hawkes processes is much more involved.
The impact of an event on the arrival rate of future events depends not only on the impact of the event itself but also on the impact of the child events and their marks.
To analyze the intensity function we therefore link marked Hawkes processes to Hawkes random measures as introduced in Horst and Xu \cite{HorstXu2019}, and then use arguments in Xu \cite{Xu2018} to give a new stochastic Volterra representation for the intensity process in terms of a martingale measure.  The stochastic Volterra integrals in this new representation can be approximated by martingales. This allows us to prove a functional CLT for the cumulative intensity by proving the weak convergence of these martingales. Specifically, we prove that 
\beqnn
\frac{1}{T} \int_0^{Tt} Z(s)ds = C t +   \frac{C}{\sqrt{T}} B(t) + o(1/\sqrt{T}),\quad \mbox{a.s. for large $T>0$}
\eeqnn
where $\{B(t): t\geq 0  \}$ is a standard Brownian motion. 
In the second step, we use the functional CLT for the cumulative intensity to prove the weak convergence of the normalized marked Hawkes point measure with immigration to a measure-valued process. 
The limit process is given in terms of a Gaussian white noise and a lifting map of a one-dimensional Brownian motion associated with the probability measure $\nu_H (du)$. Specifically, for $T>0$ large enough,
\beqnn
\frac{1}{T}N_H(dTt,du)\sim C dt\nu_H(du)+ \frac{C}{\sqrt{T}} W(dt,du) +  \frac{C}{\sqrt{T}} dB(t)\cdot \nu_H(du),
\eeqnn
where $W(dt,du)$ is a Gaussian white noise and $\{B(t): t\geq 0  \}$ is the Brownian motion from the approximation of the cumulative intensity process. 
The Gaussian white noise and the Brownian motion a correlated, due to the self-exciting property of event arrivals. 

Our second main contribution is to provide functional limit theorems for the associated shot noise process that describes the impact of the events. 
We assume that the shot shape function of the $j$-th event in $N$ (resp. $N'$) is $\psi(\cdot-\tau_j,\xi_j)$ (resp. $\psi(\cdot-\sigma_j,\eta_j)$), where $\psi: \mathbb{R}_+\times\mathbb{U}\mapsto \mathbb{R}$ is right continuous with left limits in $t$ and define the corresponding shot noise process by
\beqlb\label{eqn1.04}
S^\psi(t) \ar:=\ar  \sum_{j=1}^{N_t} \psi(t-\sigma_j,\eta_j) + \sum_{k=1}^{N'_t} \psi (t-\sigma_k,\xi_k).
\eeqlb
The second term on the right side of the above equality is a Poisson shot noise process, which has been widely applied to e.g.~bunching in traffic \cite{Bartlett1963}, computer failure times \cite{Lewis1964}, earthquake aftershocks \cite{Vere-Jones1970}, insurance \cite{KluppelbergMikosch1995a}, finance \cite{Samorodnitsky1996} and workload input models  \cite{MaulikResnick2003}. 
When the correlation function $\mathbf{E}[\psi (kt, \eta)\psi(ks,\eta)]$ is regularly varying as $k\to \infty$ for any $t,s\geq 0$, Kl\"uppelberg and Mikosch \cite{KluppelbergMikosch1995b} proved the  weak convergence of the normalized Poisson shot noise process to a self-similar Gaussian process, which is  a Brownian motion when the shot shape function is light-tailed, i.e., ${\mathbf E}[\psi (t,\eta)]= C+ o(1/\sqrt{t})$. 
In this paper, under a light-tailed condition for the shot shape function, we prove that the normalized Hawkes shot noise process with random marks and homogeneous immigration converges weakly to a Brownian martingale.


Marked Hawkes point measures are tailor made to study the dynamics of budding microbial populations with immigration in a host and their interaction with that host. In this application the mark of a microbe comprises its life length (which is rarely exponentially distributed as argued by e.g.~Holbrook and Menninger \cite{HolbrookMenninger2002} and Wood et al. \cite{Wood2004}) and the type of toxin it releases. The shot shape function describes the relaxation of the toxin as a function of the age of the microbe and the corresponding shot noise process describes the cumulative relaxation of toxin or damage made to the host by the entire population of microbes at any given point in time. Our description of budding microbial populations is consistent to Peter Jagers' \cite{Jagers2010} suggestion that biological populations should be finite and individual based. We prove a functional central limit theorem for the toxin cumulative process. When the microbes release toxins at a unit rate, then the toxin cumulative process reduces to the integral of microbial population. Pakes \cite{Pakes1975} proved a central limit theorem for the integral of microbial population; we obtain a corresponding {\it functional} central limit theorem. 

The remainder of this paper is organized as follows. In Section 2 provides an integral representation of marked Hawkes point measure and its intensity process. The main results are given in Section 3. The proof of the functional central limit theorems is given in Section 4. The proofs of the functional CLTs do not use the corresponding LLNs; the functional LLNs turn out to be immediate corollaries of the functional CLTs. The application of marked Hawkes point measures to budding microbial populations in a host is given in Section 5.

\textbf{Notation.} 
For any functions $F,G$ on $\mathbb{R}$, denote by $F*G$ the convolution of $F$ and $G$, and $F^{(*n)}$ the $n$-th convolution of $F$.  We make the convention that for any $t_1\leq t_2\in\mathbb{R}$
\beqnn
\int_{t_1}^{t_2}=\int_{(t_1,t_2]}\quad \mbox{and} \quad \int_{t_1}^{\infty}=\int_{(t_1,\infty)}.
\eeqnn
For any functions $f(t)$ on $\mathbb{R}_+$ and $g(t,u)$ on $\mathbb{R}_+\times \mathbb{U}$, let  $ \|f\|_\infty:= \sup_{t\geq 0}|f(t)| $, $  \|g(u)\|_\infty:= \sup_{t\geq 0}|g(t,u)|$  and 
\beqnn 
\|f\|_{L^{\kappa}}:= \int_0^\infty |f(t)|^\kappa dt
\quad \mbox{and}\quad 
\|g(u)\|^\kappa_{L^{\kappa}}:= \int_0^\infty |g(t,u)|^\kappa dt,\quad \kappa>0.
\eeqnn
We denote by $B(\mathbb{U})$ be the space of bounded Borel functions on $\mathbb{U}$ and $C_b(\mathbb{U})$ be the subspace of continuous elements of $B(\mathbb{U})$.
 Let $\mathcal{M}(\mathbb{U})$ be the space of finite Borel measures on $\mathbb{U}$ endowed with the weak convergence topology, i.e. for $\{ \nu_n \}_{n\geq 1}$,  $\nu\in \mathcal{M}(\mathbb{U})$, we say that $\nu_n\to \nu$ if $\nu_n(f)\to \nu(f)$ for any $f\in C_b(\mathbb{U})$, where $\nu(f)$ denotes the integral of a function $f$ with respect to a measure $\nu$ if the integral exists. 
 Let $\mathcal{S}(\mathbb{U})$ be the space of finite Borel signed measures on $\mathbb{U}$, which is also endowed with the weak convergence topology.


   \section{Stochastic integral representations}
  \setcounter{equation}{0}
  \medskip

 In this section, we give a stochastic integral representations for marked Hawkes point measures and their intensity processes that will be important for the subsequent analysis of our functional limit theorems. 
 From (\ref{eqn1.02}) and the independence of $\xi_k$ and $\{\tau_i:i=1,\cdots, k \}$ for any $k\geq1$, we see that the random point measure $N_H(ds,du)$ defined by (\ref{eqn2.03}) has the intensity $Z(s-)ds\nu_{H}(du)$. That is, for any $f\in B(\mathbb{U})$, 
\beqlb\label{eqn1.01}
\mathbf{E}\Big[\int_0^t \int_{\mathbb{U}} f(u) N_H(ds,du)  \Big]\ar=\ar  \mathbf{E}\Big[  \int_0^t \int_{\mathbb{U}} f(u)Z(s)ds \nu_H(du) \Big],\quad t\geq 0.
\eeqlb

 We denote by $\{N_{I,t}(A):t\geq 0,A\in\mathscr{U}\}$ and by $\{N_{H,t}(A): t\geq 0,A\in\mathscr{U}\}$ the $\mathcal{M}(\mathbb{U})$-valued processes associated to $N_I(ds,du)$ and $N_H(ds,du)$, respectively.
Following the argument in \cite[p.93]{IkedaWatanabe1989}, on an extension of the original probability space we can define a time-homogeneous Poisson random measure $N_0(ds,du,dz)$  on $(0,\infty)\times \mathbb{U}\times \mathbb{R}_+$ with intensity $ds\nu_H(du)dz$ such that $N_0(ds,du,dz)$ is independent of $N_I(ds,du)$ and  
\beqlb\label{eqn2.04}
N_{H,t}(f)\ar=\ar \int_0^t \int_{\mathbb{U}}\int_0^{Z(s-)} f(u)N_0(ds,du,dz),\quad f\in B(\mathbb{U}).
\eeqlb
 We can thus rewrite the intensity process $\{Z(t):t\geq 0 \}$ defined by (\ref{eqn1.02})-(\ref{eqn1.03}) as follows: for any $t\geq 0$,
 \beqlb\label{eqn2.05}
 Z(t) \ar=\ar \mu_0(t)+ \int_0^t \int_\mathbb{U} \phi(t-s, u)N_I(ds,du) +\int_0^t \int_\mathbb{U}\int_0^{Z(s-)} \phi(t-s, u)N_0(ds,du,dz) . 
 \eeqlb
 We assume throughout  that the functions $\phi_H(t) := \mathbf{E}[\phi(t,\xi_1)]$ and $\phi_I(t) :=\mathbf{E}[\phi(t,\eta_1)]$ are integrable on $[0,\infty)$.
 Taking expectations on the both sides of (\ref{eqn2.05}), we have
 \beqlb\label{eqn2.07}
 \mathbf{E}[Z(t)]\ar=\ar \mathbf{E}[\mu_0(t)] + \lambda_I \int_0^t \phi_I(s)ds + \int_0^t \phi_H(t-s)\mathbf{E}[Z(s)]ds,
 \eeqlb
 which is a linear Volterra integral equation. 
 By Theorem~3.5 in \cite[p.44]{Gripenberg1990}, the unique solution is given by
 \beqlb\label{eqn2.08}
 \mathbf{E}[Z(t)]
 \ar=\ar \mathbf{E}[\mu_0(t)] +\int_0^t R_H(t-s)\mathbf{E}[\mu_0(s)]ds  + \lambda_I \int_0^t R_I(s)ds,
 \eeqlb
 where $R_H(\cdot)$ is the resolvent kernel associated with  $\phi_H(\cdot)$ defined by the Volterra integral equation
 \beqlb\label{eqn2.09}
 R_H(t)\ar=\ar \phi_H(t) + \int_0^t R_H(t-s)\phi_H(s)ds= \phi_H(t) + \int_0^t \phi_H(t-s)R_H(s)ds
 \eeqlb
 and
 \beqlb\label{eqn2.10}
 R_I(t)\ar=\ar \phi_I(t) + \int_0^t R_H(t-s) \phi_I(s)ds .
 \eeqlb
Integrating the both sides of (\ref{eqn2.09}) and (\ref{eqn2.10}) over the interval $(0,\infty)$, we have
 \beqlb\label{eqn2.11}
 \|R_H\|_{L^1}=\|\phi_H\|_{L^1}\cdot (1+  \| R_H\|_{L^1} ) 
 \quad\mbox{and}\quad
  \| R_I\|_{L^1}=  \|\phi_I\|_{L^1} \cdot (1+  \| R_H\|_{L^1} ),
 \eeqlb
 from which we see that $ \|R_H\|_{L^1}+\|R_I\|_{L^1}<\infty$ if and only if $\|\phi_H\|_{L^1} <1$. In this case,
  \beqlb\label{eqn2.12}
 \|R_H\|_{L^1}=\frac{\|\phi_H\|_{L^1}}{1-\|\phi_H\|_{L^1}}<\infty\quad \mbox{and}\quad \| R_I\|_{L^1}= \frac{\|\phi_I\|_{L^1}}{1-\|\phi_H\|_{L^1}} <\infty.
 \eeqlb
 
 \begin{lemma}\label{Thm201}
  If $\sup_{s\geq 0}\mathbf{E}[\mu_0(s)] + \|\phi_I\|_{L^1}<\infty$, then $\mathbf{E}[Z(\cdot)]$ is uniformly bounded if and only if $ \|\phi_H\|_{L^1}<1$.
 	Moreover, in this case,
 	\beqlb\label{eqn2.13}
 	\sup_{t\geq 0}\mathbf{E}[Z(t)]\leq \frac{\sup_{s\geq 0} \mathbf{E}[\mu_0(s)]+\lambda_I\cdot \|\phi_I\|_{L^1} }{1-\|\phi_H\|_{L^1}}.
 	\eeqlb
 \end{lemma}
 \proof From (\ref{eqn2.07}), the uniform boundedness of $\mathbf{E}[Z(t)]$ induces immediately that $\| R_H\|_{L^1} + \| R_I\|_{L^1}<\infty$ and $ \|\phi_H\|_{L^1}<1$.
 For the converse, from  (\ref{eqn2.07}) and (\ref{eqn2.12}),
 \beqlb\label{eqn2.14}
 \sup_{t\geq 0}\mathbf{E}[Z(t)] \leq \sup_{t\geq 0} \mathbf{E}[\mu(t)] \big(1+ \|R_H\|_{L^1}\big) +  \lambda_I\cdot \|R_I\|_{L^1}
 = \frac{\sup_{t\geq 0} \mathbf{E}[\mu(t)]+\lambda_I\cdot \|\phi_I\|_{L^1} }{1-\|\phi_H\|_{L^1}}< \infty.
 \eeqlb
 \qed

 
 \section{Main results}
  \setcounter{equation}{0}
 \medskip
 
 In this section we state our functional limit theorems for the $\mathcal{M}(\mathbb{U})$-valued processes $\{N_{H,t}(A) :t\geq 0, A\in\mathscr{U}  \}$ and $\{N_{I,t}(A) :t\geq 0, A\in\mathscr{U}  \}$ and their shot noise processes  under the stability condition $\|\phi_H\|_{L^1}<1$. 
%
%
 
 \subsection{Functional laws of large numbers}

%
 
 \subsubsection{Point measures}
 For any $T> 0$ and $i\in\{H,I \}$, we define the rescaled measure-valued process $\{N^T_{i,t}(A) :t\geq 0, A\in\mathscr{U}  \}$  by
 \beqnn
 N^T_{i,t}(A):=  \frac{1}{T}N_{i,Tt}(A).
 \eeqnn
 
 The asymptotic analysis of the Poisson random measure  $N_I(ds,du)$ is  standard; see Theorem~7.10 in Walsh  \cite{Walsh1986}.
 Let  $\tilde{N}_I(ds,du):=N_I(ds,du)-\lambda_Ids\nu_I(du)$  be the compensated point measure of $N_I(ds,du)$. 
 From the Burkholder-Davis-Gundy inequality, for any $f\in  B(\mathbb{U})$,  
 \beqlb\label{eqn3.01}
 \mathbf{E}\Big[\sup_{t\in[0,1]} \Big|\int_0^{Tt}\int_{\mathbb{U}}\frac{f(u)}{T}\tilde{N}_I(ds,du) \Big|^2\Big]\leq   \frac{C}{T}\int_{\mathbb{U}} |f(u)|^2\nu_I(du) \to 0, \quad \mbox{as } T\to\infty,
 \eeqlb
 which immediately yields the following functional law of large numbers for $\{N_{I,t}(A):t\geq 0, A\in\mathscr{U}  \}$.
 \begin{lemma}\label{Thm301}
As $T\to\infty$, the recaled $\mathcal{M}(\mathbb{U})$-valued process $\{N^T_{I,t}(A):t\geq 0,A\in\mathscr{U}\}$ converges to the deterministic $\mathcal{M}(\mathbb{U})$-valued process $ \{ \lambda_I\cdot t\cdot \nu_I(A):t\geq 0,A\in\mathscr{U}\}$ uniformly in probability on any bounded time interval.
 \end{lemma}

 We now consider the asymptotic behavior of  the rescaled Hawkes random measure $\{N^T_{H,t}(A):t\geq 0,A\in\mathscr{U}\}$. 
 Taking expectations on the both sides of (\ref{eqn2.04}), we have 
 \beqlb\label{eqn3.04}
 \mathbf{E}\Big[ \int_0^t \int_{\mathbb{U}}f(u)N_H(ds,du) \Big]\ar=\ar   \nu_H(f) \int_0^t\mathbf{E}[Z(s)] ds .
 \eeqlb
 From (\ref{eqn2.08}) and Fubini's lemma,
 \beqlb\label{eqn3.05}
 \mathbf{E}\Big[\int_0^t Z(s)ds  \Big]=\int_0^t \mathbf{E}[\mu_0(s)]ds +\int_0^t R_H(t-s)ds\int_0^s \mathbf{E}[\mu_0(r)]dr + \int_0^t \lambda_I ds \int_0^sR_I(r)dr.
 \eeqlb
 When $\mathbf{E}[\mu_0(\cdot)]$ is integrable on $[0,\infty)$, we may conjecture that   $\mathbf{E}[\int_0^{t}Z(s)ds]
 \approx \lambda_I\cdot \|R_I\|_{L^1}\cdot t$ for $t>0$ large enough, because as $t\to \infty$,
 \beqnn
  \mathbf{E}\Big[\int_0^t Z(s)ds  \Big] \sim \lambda_I \int_0^t (t-s)R_I(s)ds + o(1)\sim \lambda_I\cdot \|R_I\|_{L^1} \cdot t + o(1). 
 \eeqnn
 To obtain the exact rate of convergence, we need the following moment assumption on $\mu_0(\cdot)$ and  $\phi(\cdot)$.
 \begin{condition} \label{C1}
 	There exist constants $\alpha> 1$ and $\theta_0> \frac{\alpha}{2\alpha-2}$ such that for $i \in \{H,I \}$,
 	\beqlb\label{eqn3.06}
 	\sup_{t\geq 0}	\mathbf{E}[|\mu_0(t)|^{2\alpha}] + \mathbf{E}[\|\mu_0\|_{L^{2\alpha}}^{2\alpha}]  <\infty
 	\eeqlb
 	and
 	\beqlb\label{eqn3.07}
 	\int_0^\infty t^{\theta_0}\phi_i(t)dt+ \int_{\mathbb{U}}\left(\|\phi(u)\|_\infty^{2\alpha}+ \|\phi(u)\|_{L^1}^{2\alpha}\right)\nu_i(du)<\infty.
 	\eeqlb
 \end{condition}
 
 Under the preceding condition,  $\sup_{t\geq 0}\mathbf{E}[|Z(t)|^{2\alpha}]<\infty$; see  Corollary~\ref{ThmA04}. 
We shall assume $\alpha\in(1,2]$; the case $\alpha>2$ is much simpler to analyze.
 From (\ref{eqn3.05}), 
 \beqlb\label{eqn3.08}
 \mathbf{E}\Big[\int_0^{t}Z(s)ds \Big] - \lambda_I\cdot \|R_I\|_{L^1}\cdot t
 \ar=\ar \int_0^t \mathbf{E}[\mu_0(s)]ds-   \int_0^{t} \lambda_I ds \int_s^\infty R_I(r)dr\cr
 \ar\ar   +\int_0^t R_H(t-s)ds\int_0^s \mathbf{E}[\mu_0(r)]dr   .
 \eeqlb
 From (\ref{eqn2.12}), we have for any $T\geq 0$,
 \beqlb\label{eqn3.09}
  \sup_{t\in[0,T]}\Big|\mathbf{E}\Big[\int_0^{t}Z(s)ds \Big]
 - \lambda_I\cdot \|R_I\|_{L^1}\cdot t\Big|
 \ar\leq\ar  \mathbf{E}[\|\mu_0\|_{L^1} ](1+\|R_H\|_{L^1} )+ \int_0^{T}\lambda_I dt \int_t^\infty R_I(s)ds\cr 
 \ar\leq\ar \frac{ \mathbf{E}[\|\mu_0\|_{L^1} ] }{1-  \|\phi_H\|_{L^1}} + \int_0^{T}\lambda_I dt \int_t^\infty R_I(s)ds. 
 \eeqlb
  Changing the order of integration in the second term on the right side of the last inequality, for any $\kappa\in(0,\theta_0\wedge 1)$,  
 \beqlb\label{eqn3.09.01}
 \int_0^{T}dt \int_t^\infty R_I(s)ds\ar=\ar  \int_0^{T} s R_I(s)ds + T\int_T^\infty R_I(s)ds\cr
 \ar\leq\ar  T^{1-\kappa}\int_0^{T} s^\kappa R_I(s)ds + T^{1-\kappa}\int_T^\infty s^\kappa R_I(s)ds\cr
 \ar=\ar T^{1-\kappa}\int_0^\infty s^\kappa R_I(s)ds.
 \eeqlb
 Moreover,  from (\ref{eqn2.09}) and the inequality $|a+b|^\kappa\leq a^\kappa+b^\kappa$ for any $a,b\geq 0$,
 \beqlb\label{eqn3.10}
\int_0^\infty s^\kappa R_H(s)ds\ar=\ar    \int_0^\infty s^\kappa \phi_H(s)ds + \int_0^\infty s^\kappa ds\int_0^s R_H(s-r)\phi_H(r)dr   \cr
 \ar\leq \ar   \int_0^\infty s^\kappa\phi_H(s)ds   +  \int_0^\infty ds  \int_0^s R_H(s-r)\cdot r^\kappa\phi_H(r)dr \cr 
  \ar\ar  +   \int_0^\infty ds\int_0^s (s-r)^\kappa R_H(s-r)\cdot \phi_H(r)dr   \cr
 \ar=\ar  \int_0^\infty r^\kappa \phi_H (r)dr \cdot\big[ 1 + \| R_H\|_{L^1} \big] +  \|\phi_H\|_{L^1}\cdot \int_0^\infty s^\kappa R_H(s)ds.
 \eeqlb
 Solving this inequality, we conclude from (\ref{eqn3.07})  that
 \beqlb\label{eqn3.11}
 \int_0^\infty s^\kappa R_H(s)ds\ar\leq\ar \frac{\int_0^\infty r^\kappa \phi_H (r)dr  }{|1- \|\phi_H\|_{L^1}|^2} <\infty.
 \eeqlb
 Similarly, we also have
 \beqlb\label{eqn3.12}
 \int_0^\infty s^\kappa R_I(s)ds\ar\leq \ar   \int_0^\infty r^\kappa\phi_I (r)dr \cdot\big[ 1 + \| R_H\|_{L^1} \big] +  \|\phi_I\|_{L^1}\cdot \int_0^\infty s^\kappa R_H(s)ds \cr
 \ar\ar\cr
 \ar\leq\ar   \frac{ \int_0^\infty r^\kappa\phi_I (r)dr }{ 1- \|\phi_H\|_{L^1} }  +  \frac{ \|\phi_I\|_{L^1}\cdot \int_0^\infty r^\kappa\phi_H (r)dr}{|1- \|\phi_H\|_{L^1}|^2}  <\infty.
 \eeqlb
 Taking this and (\ref{eqn3.09.01}) back into (\ref{eqn3.09}), we get
 \beqnn
 	\sup_{t\in[0,T]}\Big|\mathbf{E}\Big[\int_0^{t}Z(s)ds \Big]
 -\lambda_I\cdot \|R_I\|_{L^1}\cdot t \Big| \leq C(1+T^{1-\kappa}),
 \eeqnn
 which induces the following proposition.
 \begin{proposition}\label{Thm302}
 	Under Condition~\ref{C1}, we have 
 	\beqlb\label{eqn3.14}
 	\lim_{T\to\infty}\sup_{t\in[0,1]}\Big| \mathbf{E}\Big[\frac{1}{T}\int_0^{Tt} Z(s)ds  \Big] -   \lambda_I\cdot \|R_I\|_{L^1}\cdot t \Big|= 0.
 	\eeqlb
 \end{proposition}
 
 
 The previous proposition shows that the law of large numbers holds for the cumulative intensity $\int_0^{\cdot} Z(s)ds$. 
The following lemma estibalishes a \textit{functional} law of large numbers for this process. 
The proof is given in Section~\ref{ProofA}. 
  \begin{lemma} \label{Thm303}
 	Under Condition~\ref{C1}, we have  as $T\to\infty$,
 	\beqlb\label{eqn3.15}
 	\sup_{t\in[0,1]} \Big| \frac{1}{T}\int_0^{Tt} Z(s)ds - \lambda_I\cdot \|R_I\|_{L^1}\cdot t \Big|\overset{\mathbf{P}}\longrightarrow 0.
 	\eeqlb
 \end{lemma}
 
 Let us now turn to the functional law of large numbers for the marked Hawkes point measure $\{ N_{H,t}(A):t\geq 0,A\in\mathscr{U} \}$.  
Let $\tilde{N}_H(ds,du):= N_H(ds,du)-Z(s-)ds\nu_H(du)$ be the compensated point measure of $N_H(ds,du)$. For any $f\in B(\mathbb{U})$, define
 \beqnn
   \tilde{N}^T_{H,t}(f):=    \int_0^t \int_{\mathbb{U}}\frac{f(u)}{T}\tilde{N}_H(dTs,du),
 \eeqnn
 which is a martingale. Moreover, 
 \beqlb\label{eqn3.23.01}
  N^T_{H,t}(f)- \lambda_I \cdot \|R_I\|_{L^1}\cdot t\cdot \nu_H(f) \ar=\ar \tilde{N}^T_{H,t}(f)+   \Big[  \int_0^{t} Z(Ts)ds   -\lambda_I \cdot \|R_I\|_{L^1}\cdot t \Big]\cdot \nu_H(f).
 \eeqlb
 Applying the Burkholder-Davis-Gundy inequality to $\{  \tilde{N}^T_{H,t}(f): t\geq 0 \}$, from Lemma~\ref{Thm201} we have 
 \beqnn
 \mathbf{E}\Big[\sup_{t\in[0,1]}|\tilde{N}^T_{H,t}(f)|^2 \Big]\leq C  \int_0^1 \mathbf{E}[Z(Ts)]ds \int_{\mathbb{U}}\frac{|f(u)|^2}{T}\nu_H(du)\leq \frac{C}{T},
 \eeqnn
 which vanishes as $T\to\infty$. From this and Lemma~\ref{Thm303}, we can get the following functional law of large numbers for the $\mathcal{M}(\mathbb{U})$-valued process  $\{ N_{H,t}(A):t\geq 0, A\in\mathscr{U} \}$. 
 \begin{theorem}\label{Thm304}
 	Under Condition~\ref{C1},  the rescaled $\mathcal{M}(\mathbb{U})$-valued process $\{N^T_{H,t}(A):t\geq 0,A\in\mathscr{U}\}$ converges to the deterministic $\mathcal{M}(\mathbb{U})$-valued process $ \{\lambda_I\cdot \|R_I\|_{L^1}\cdot t\cdot \nu_H(A):t\geq 0,A\in\mathscr{U}\}$ uniformly in probability on any bounded time interval as $T\to\infty$. 
 \end{theorem}

 \subsubsection{Shot noise processes}
 
 Before giving the limit theorems for the shot noise processes driven by the  random point measures $N_H(ds,du)$ and $N_I(ds,du)$, we introduce some conditions on the shot shape functions.
 For any mark $u\in\mathbb{U}$, we may always assume that the total cumulative impact $\psi(\infty,u):=\lim_{t\to\infty}\psi(t,u)$ is finite. 
 Moreover, conditioned on shot shape function $\psi(\cdot)$, we also assume that both the  mean impact of an event up to age $t$ and its total mean impact are finite, i.e., for any $ t\in[0,\infty]$ and $i\in\{ H, I \}$,
 \beqlb \label{eqn3.12.01}
 \psi_i(t):=\int_\mathbb{U}\psi(t,u)\nu_i(du)<\infty .
 \eeqlb
 For $u\in\mathbb{U}$ and $t\geq 0$, we define the following functions that represent  the total impact of an event after age $t$:
 \beqlb\label{eqn3.12.02}
 \psi^{\rm c}(t,u):=\psi(\infty,u)-\psi(t,u)\quad \mbox{and}\quad \psi^{\mathrm{c}}_i(t):= \psi_i(\infty)-\psi_i(t),
 \eeqlb

 \begin{condition}\label{C2}
  There exists a constant $\theta_1>\frac{2\alpha-1}{2\alpha-2}$ such that for $i\in  \{H,I \}$,
 	\beqlb\label{eqn3.17}
 	\sup_{t\geq 0}\int_{\mathbb{U}}\Big[|\psi(t,u)|^{2\alpha} +  t^{\theta_1}\cdot \sup_{s\geq t}|\psi_i^{\mathrm{c}}(s,u)| \Big]\nu_i(du)<\infty.
 	\eeqlb
 \end{condition}

 Denote by $\{ S_H^\psi(t):t\geq 0 \}$ and $\{ S_I^\psi(t):t\geq 0 \}$  the two shot noise processes on the right side of (\ref{eqn1.04}). 
 From (\ref{eqn2.03}) and (\ref{eqn2.01}), we derive the following stochastic integral representations:
 	\beqlb\label{eqn3.18}
 	S_i^\psi(t) \ar=\ar   \int_0^t \int_\mathbb{U}  \psi(t- s, u)N_i(ds,du),\quad i\in\{H,I \}.
 	\eeqlb
 Taking expectations on the both sides of this equation with $i=I$, we have 
 \beqlb\label{eqn3.19}
  \mathbf{E}[S_I^\psi(t) ]= \lambda_I\cdot \int_0^t  \psi_I(s)ds 
  \quad \mbox{and}\quad   
  \mathbf{E}[S_I^\psi(t) ]-\psi_I(\infty)\cdot \lambda_I\cdot  t =  - \lambda_I\cdot \int_0^t  \psi^{\rm c}_I(s)ds, 
 \eeqlb
 which is uniformly bounded.
 Hence the following limit holds uniformly:
 \beqlb\label{eqn3.20}
 \lim_{T\to\infty}\mathbf{E}\big[S_I^\psi(Tt)/T \big] \ar=\ar \psi_I(\infty)\cdot \lambda_I \cdot t.
 \eeqlb
 Taking expectations on the both sides of (\ref{eqn3.18}) with $i=H$,  from  (\ref{eqn2.08})  we have 
 \beqlb\label{eqn3.21}
 \mathbf{E}[S_H^\psi(t) ] 
 \ar=\ar    \int_0^t R_H(t-s)ds \int_0^s \psi_H(s-r)\mathbf{E}[\mu_0(r)]dr \cr
 \ar\ar + \int_0^t \psi_H(t-s)\mathbf{E}[\mu_0(s)] ds  + \lambda_I \cdot\int_0^t \psi_H(t-s)ds\int_0^s R_I(r)dr.
 \eeqlb
 From Condition~\ref{C1} and \ref{C2}, we can see that the first two terms on the right side of the equality above can be uniformly bounded by 
 \beqlb
 \|\psi_H\|_\infty \cdot\int_0^\infty \mathbf{E}[\mu_0(s)]ds \cdot \big( \|R_H\|_{L^1}+1 \big)<\infty.
 \eeqlb
 Moreover,  from Condition~\ref{C1}, \ref{C2} and (\ref{eqn3.11})-(\ref{eqn3.12}), 
 \beqnn
   \big|\mathbf{E}[S_H^\psi(t) ]-  \psi_H(\infty)\cdot\lambda_I\cdot \|R_I\|_{L^1}\cdot t \big|
\ar\leq \ar
 C  + \lambda_I\|R_I\|_{L^1} \int_0^\infty |\psi_H^{\mathrm{c}}(s)| ds  +  \lambda_I  \|\psi_H\|_\infty  \int_0^t ds\int_s^\infty  R_I(r)dr. 
\eeqnn
 From (\ref{eqn3.17}), we see that the second term on the right side of the inequality above is finite. From (\ref{eqn3.09.01}) and (\ref{eqn3.12}),  we have for any $T>0$ and $\kappa\in (0,\theta_0\wedge 1 )$
\beqnn
\sup_{t\in[0,T]}\big|\mathbf{E}[S_H^\psi(t) ]-  \psi_H(\infty)\cdot\lambda_I\cdot \|R_I\|_{L^1}\cdot t \big|
\ar\leq \ar
C (1 + T^{1-\kappa}),
\eeqnn
which induces the following  convergence:
 \beqlb\label{eqn3.22}
\lim_{T\to\infty}\sup_{t\in[0,1]}\Big| \mathbf{E}\Big[\frac{1}{T}S_H^\psi(Tt)\Big]- \psi_H(\infty)\cdot \lambda_I\cdot \|R_I\|_{L^1} \cdot t\Big| =0.
\eeqlb
 From (\ref{eqn3.20}) and (\ref{eqn3.22}), we derive the following functional laws of large numbers for the shot noise processes. 
 
 \begin{theorem}\label{Thm305}
  Under Condition~\ref{C1} and \ref{C2},  we have  as $T\to\infty$,
	\beqlb\label{eqn3.23}
 	\sup_{t\in [0,1]} \Big|  \frac{1}{T}S_I^\psi (Tt)-   \psi_I(\infty)\cdot \lambda_I \cdot  t\Big| + \sup_{t\in [0,1]} \Big|  \frac{1}{T}S_H^\psi (Tt)-   \psi_H(\infty)\cdot \lambda_I\cdot  \|R_I\|_{L^1}\cdot t\Big| 		\overset{\mathbf{P}}\longrightarrow 0.
	\eeqlb
 \end{theorem}

 \subsection{Functional central limit theorems}
 
%
 \subsubsection{Point measures}
 
Lemma~\ref{Thm301} and Theorem~\ref{Thm304} show that for $T> 0$ large enough the rescaled $\mathcal{M}(\mathbb{U})$-valued processes $\{N^T_{I,t}(A):t\geq 0, A\in\mathscr{U}\}$ and $\{N^T_{H,t}(A):t\geq 0, A\in\mathscr{U}\}$ can be approximated by the deterministic processes $\{\lambda_I\cdot t\cdot \nu_I(A):t\geq 0, A\in\mathscr{U} \}$ and $\{\lambda_I\cdot \|R_I\|_{L^1}\cdot t \cdot \nu_H(A):t\geq 0, A\in\mathscr{U} \}$ uniformly in probability, respectively. Let us denote by $\{\bar{N}^T_{I,t}(A):t\geq 0, A\in\mathscr{U} \}$ and $\{\bar{N}^T_{H,t}(A):t\geq 0, A\in\mathscr{U} \}$  the respective error processes
 \beqlb\label{eqn3.24}
 \bar{N}^T_{I,t}(A):= N^T_{I,t}(A)-\lambda_I\cdot t\cdot \nu_I(A)
  \quad \mbox{and}\quad
 \bar{N}^T_{H,t}(A):= N^T_{H,t}(A)- \lambda_I\cdot \|R_I\|_{L^1}\cdot t\cdot \nu_H(A).
 \eeqlb
 From (\ref{eqn2.01}), we see that  
 \beqnn
 \sqrt{T}\bar{N}^T_{I,t}(A)=\int_0^t\int_{A} \frac{1}{\sqrt{T}}\tilde{N}_I(dTs,du),
 \eeqnn
 which is a worthy martingale measure in the sense of Walsh \cite[p.291]{Walsh1986}.
 The following lemma gives a functional central limit theorem for $\{N_{I,t}(A):t\geq 0, A\in\mathscr{U} \}$. 
 The proof will be given in Section~\ref{ProofA}.
 
 \begin{lemma}\label{Thm401}
 	As $T\to\infty$, we have that $\{ \sqrt{T}\bar{N}^T_{I,t}(A):t\geq 0, A\in\mathscr{U}  \}$ converges to $\{W_{I,t}(A): t\geq 0, A\in\mathscr{U} \}$
 	weakly in the space $\mathbb{D}([0,\infty),\mathcal{S}(\mathbb{U}))$, where $\{W_{I,t}(A): t\geq 0, A\in\mathscr{U} \}$ is a Gaussian white noise on  $\mathbb{U}$ with intensity $\lambda_I\cdot dt\nu_I(du)$.
 	
 \end{lemma}
 
 Next, we consider the functional central limit theorem for $\{N_{H,t}(A):t\geq 0, A\in\mathscr{U} \}$.  
 Unlike  $\{ \sqrt{T}\bar{N}^T_{I,t}(A):t\geq 0, A\in\mathscr{U}  \}$, the corresponding measure $\{ \sqrt{T}\bar{N}^T_{H,t}(A):t\geq 0, A\in\mathscr{U}  \}$ is not a martingale measure. 
 Loosely speaking, the limit process is a sum of a Gaussian white noise and a lifting map of a Brownian motion associated with the probability measure $\nu_H(du)$ resulting from the cumulative intensity. 
 Specifically, from (\ref{eqn3.24}) and (\ref{eqn3.23.01}),  
 \beqlb\label{eqn4.03}
 \bar{N}^T_{H,t}(A)\ar=\ar    \tilde{N}^T_{H,t}(A)
 +  \Big[  \int_0^{t} Z(Ts)ds   -\lambda_I \cdot \|R_I\|_{L^1}\cdot t \Big]\cdot \nu_H(A).
 \eeqlb
 The following proposition shows the  weak convergence of $\{\tilde{N}^T_{H,t}(A):t\geq 0, A\in\mathscr{U} \}$ to a Gaussian white noise as $T\to\infty$; a detailed proof will be given in Section~\ref{ProofA}.
 \begin{proposition}\label{Thm402}
 	Under Condition~\ref{C1}, as $T\to\infty$ we have that $\{ \sqrt{T}\tilde{N}^T_{H,t}(A):t\geq 0, A\in\mathscr{U}   \}$ converges weakly to $\{  W_{H,t}(A): t\geq0, A \in \mathscr{U}  \}$ in the space $\mathbb{D}([0,\infty),\mathcal{S}(\mathbb{U}))$,
  where $\{  W_{H,t}(A): t\geq0, A \in \mathscr{U}  \}$  is a Gaussian white noise on  $\mathbb{U}$ with intensity $\lambda_I \cdot \|R_I\|_{L^1}\cdot dt\nu_H(du)$ and independent of $\{  W_{I,t}(A): t\geq0, A \in \mathscr{U}  \}$.
 	
 \end{proposition}
 
 Before giving the key result about the weak convergence of the second term on the right side of (\ref{eqn4.03}), we introduce the two-parameter function 
 \beqlb\label{eqn4.05}
 R(t,u)\ar:=\ar \phi(t,u)+\int_0^t R_H(t-s) \phi(s,u)ds.
 \eeqlb
 The function can be interpreted as describing the mean impact up to time $t$ of an event with mark $u$ on the future intensity.
 For $i\in\{H,I \}$,  integrating both sides of (\ref{eqn4.05}) with respect to $\nu_i(du)$, we have $R_i(t)= \int_\mathbb{U}R(t,u)\nu_i(du)$.
 Thus $\{R(t,u):u\in\mathbb{U}\} $ can be considered as the decomposition of $R_i(t)$ on the space $\mathbb{U}$.
 From Condition~\ref{C1}, 
 \beqlb\label{eqn4.06}
 \|R(u)\|_\infty \leq \|\phi(u)\|_\infty \cdot  (1+\|R_H\|_{L^1}) =\frac{\|\phi(u)\|_\infty}{1-\|\phi_H\|_{L^1}},
 \eeqlb
 for any $u\in\mathbb{U}$.
 Integrating both sides of (\ref{eqn4.05}) over the interval $(0,\infty)$, we also have for any $u\in\mathbb{U}$,
 \beqlb\label{eqn4.07}
 \|R(u)\|_{L^1}  
 \ar=\ar \|\phi(u)\|_{L^1} \cdot  (1+\|R_H\|_{L^1} )
 = \frac{\|\phi(u)\|_{L^1}}{1-\|\phi_H\|_{L^1}}.
 \eeqlb
 In view of Condition~\ref{C1}, this implies that $\nu_i(\|R(\cdot)\|_{L^1}^2)<\infty$ and hence $W_{i,t}(\|R(\cdot)\|_{L^1})$ is well defined for $i\in\{ H,I \}$. We are now ready to state the functional central limit theorem for the cumulative intensity process with proof will be given in Section~\ref{ProofA}. 
 
 \begin{proposition}\label{Thm403}
 	Under Condition~\ref{C1},  we have as $T\to\infty$,
 	\beqlb\label{eqn4.08}
 	\sqrt{T} \Big(  \int_0^{t} Z(Ts)ds 
 	-\lambda_I \cdot \|R_I\|_{L^1}\cdot t
 	\Big) \to   W_{H,t}(\| R(\cdot) \|_{L^1})+W_{I,t}(\| R(\cdot) \|_{L^1})=\sigma_Z B_Z(t),
 	\eeqlb
 	weakly in the space $\mathbb{D}([0,\infty),\mathbb{R})$, where $\{B_Z(t):t\geq \}$ is a standard Brownian motion and
 	\beqlb\label{eqn4.08.01}
 	\sigma_Z^2= \lambda_I \cdot   \frac{\| R_I\|_{L^1}\cdot  \nu_H(\|\phi(\cdot)\|_{L^1}^2)+\nu_I(\|\phi(\cdot)\|_{L^1}^2)}{|1-\|\phi_H\|_{L^1}|^2}.
 	\eeqlb
 \end{proposition}

Combining Propositions~\ref{Thm402} and \ref{Thm403} with (\ref{eqn4.03}), we get the functional CLT for $\{N_{H,t}(A):t\geq 0, A\in\mathbb{U}  \}$. 
 
 \begin{theorem} \label{Thm404}
 	Under Condition~\ref{C1},  as $T\to\infty$ we have  $\{ 	\sqrt{T}\bar{N}_{H,t}^T (A):t\geq 0, A\in\mathscr{U}  \}$ converges to $\{W_{H,t}(A) + \sigma_Z B_Z(t) \cdot \nu_H(A):t\geq 0, A\in\mathscr{U} \}$  
 	weakly in the space $\mathbb{D}([0,\infty),\mathcal{S}(\mathbb{U}))$.
 	
 \end{theorem}
 
 \subsubsection{Shot noise processes}
 
We are now going to establish the functional limit theorems for the shot noise processes $\{S_I^\psi(t):t\geq 0\}$ and $\{S_H^\psi(t):t\geq 0\}$. 
 From Condition~\ref{C2}, we can see that $\psi(Tt,u)\sim\psi(\infty,u)$ for $T\geq 0$ large enough. 
 Thus we may approximate $S_I^\psi(t)$ and $S_H^\psi(t)$ by the  semi-martingales: 
 \beqlb\label{eqn4.10}
 N_{i,t}^T( \psi(\infty, \cdot))= \int_0^t \int_\mathbb{U}  \psi(\infty, u)\frac{1}{T}N_i(dTs,du),\quad  i\in\{H,I \}.
 \eeqlb
The error processes are given by
 \beqlb\label{eqn4.11}
 \varepsilon^\psi_{T,i}(t):=  \int_0^t \int_\mathbb{U}  \psi^{\rm c}(T(t-s), u)\frac{1}{T}N_i(dTs,du).
 \eeqlb
 These vanish as $T\to\infty$, due to the following lemma whose proof will be given in Section~\ref{AppendixB}.
 \begin{lemma}\label{Thm405}
 	Under Condition~\ref{C1} and \ref{C2},  both  $\{ \sqrt{T} \varepsilon^\psi_{T,I}(t):t\geq 0\}$ and $\{ \sqrt{T} \varepsilon^\psi_{T,H}(t):t\geq 0\}$  converge weakly to $0$ in the space $\mathbb{D}([0,\infty),\mathbb{R})$ as $T\to\infty$.
 \end{lemma}

 From (\ref{eqn4.10}), Lemma~\ref{Thm401}, \ref{Thm405} and Theorem~\ref{Thm404}, we get the functional central limit theorems for $\{S_H^\psi(t) :t\geq 0  \}$ and $\{S_I^\psi(t) :t\geq 0  \}$.
 
 \begin{theorem}\label{Thm406}
 	Under Condition~\ref{C1} and \ref{C2},  we have  as $T\to\infty$,
 	\beqlb\label{eqn4.12}
 	\sqrt{T}\Big(  \frac{1}{T}S_I^\psi (Tt)-   \psi_I(\infty)\cdot \lambda_I\cdot  t\Big)\to W_{I,t} (\psi(\infty,\cdot))  
 	\eeqlb
 	and
 	\beqlb\label{eqn4.13}
 	\sqrt{T}\Big(   \frac{1}{T}S_H^\psi (Tt)-   \psi_H(\infty)\cdot \lambda_I\cdot  \|R_I\|_{L^1}\cdot t\Big)
 	\ar\to\ar  W_{H,t} (\psi(\infty,\cdot)) +\sigma_Z\cdot \psi_H(\infty) \cdot B_Z(t) ,
 	\eeqlb
 	weakly in the space $\mathbb{D}([0,\infty),\mathbb{R})$.
 \end{theorem}

 \begin{remark}{\rm
 	We emphasis that the proof of Theorem~\ref{Thm406} follows from Lemma~\ref{Thm405} and Theorem~\ref{Thm404}, which does not require the law of large numbers (Theorem~\ref{Thm305}). In other words, Theorem~\ref{Thm305} follows directly from Theorem~\ref{Thm406}, i.e., from (\ref{eqn4.12}) and (\ref{eqn4.13}),
 	\beqnn
 	\frac{1}{T}S_I^\psi (Tt)-   \psi_I(\infty)\cdot \lambda_I \cdot  t\to 0 \quad \mbox{and}\quad  \frac{1}{T}S_H^\psi (Tt)-   \psi_H(\infty) \cdot \lambda_I\cdot  \|R_I\|_{L^1}\cdot t\to 0, 
 	\eeqnn
 	weakly in the space $\mathbb{D}([0,\infty),\mathbb{R})$ and hence uniformly in probability on any bounded interval; see  \cite[p.124]{Billingsley1999}.  }
 \end{remark}
 
%

 \subsection{Examples}
 
 In this section, we illustrate how our framework can be used to derive a functional central limit theorem for standard marked Hawkes processes. 
 Let $N_\lambda(ds,du)$ be a marked Hawkes point measure  with $(\mathscr{F}_t)$-intensity 
 \beqlb\label{eqn2.05.01}
 Z_\lambda(t) \ar=\ar \lambda +\int_0^t \int_\mathbb{U}\int_0^{Z_\lambda (s-)} \phi(t-s, u)N_0(ds,du,dz)
 \eeqlb
 for some $\lambda>0$. This intensity does not satisfy Condition \ref{C1}. Instead, we now show that there exists a marked Hawkes point measure with immigration that is equavilent to $N_\lambda(ds,du)$.
 For any $t\geq 0$, let $Z(t):=Z_\lambda(t)- \lambda$, which satisfies the following equation
  \beqnn
  Z(t)
  \ar=\ar \int_0^t \int_\mathbb{U}\int_{Z(s-)}^{Z(s-)+\lambda } \phi(t-s, u)N_0(ds,du,dz)
  +\int_0^t \int_\mathbb{U}\int_0^{Z(s-)} \phi(t-s, u)N_0(ds,du,dz).
  \eeqnn
  From the orthogonality and homogeneity of $N_0(ds,du,dz)$ in space, we can see that 
 the two processes on the right side of the last equality are independent and
 \beqnn
 N_I(ds,du):= \int_{Z(s-)}^{Z(s-)+\lambda } N_0(ds,du,dz)
 \eeqnn
 is a Poisson random measure on $(0,\infty)\times \mathbb{U}$ with intensity $\lambda\cdot ds\nu_H(du)$. Thus, we can rewrite (\ref{eqn2.05.01}) as 
 \beqnn
 Z(t)\ar=\ar  \int_0^t \int_\mathbb{U} \phi(t-s, u) N_I(ds,du)+\int_0^t \int_\mathbb{U}\int_0^{Z(s-)} \phi(t-s, u)N_0(ds,du,dz)
 \eeqnn
 and (\ref{eqn2.04}) as
 \beqlb
 N_{\lambda,t}(f)\ar=\ar \int_0^t \int_{\mathbb{U}} f(u)N_I(ds,du) +  \int_0^t \int_{\mathbb{U}}\int_0^{Z_\lambda (s-)} f(u)N_0(ds,du,dz).
 \eeqlb
 Moreover, the random point measure $N_H(ds,du):= N_0(ds,du,[0,Z(s-)))$ is a Marked Hawkes point measure with homogeneous immigration on $(0,\infty)\times\mathbb{U}$ as defined in Section~1 with $\lambda_I=\lambda$ and $\nu_I(du)=\nu_H(du)$. 
 Applying Lemma~\ref{Thm401} and Theorem~\ref{Thm404}, under Condition~\ref{C1} with $\mu_0(t)\equiv 0$ we can prove the weak convergence of the normalized standard marked Hawkes point measure defined as
 \beqnn
 \sqrt{T}\bar{N}_{\lambda,t}^T (A) := \sqrt{T} \Big(  \frac{1}{T}N_\lambda([0,Tt], A)- \frac{\lambda}{1- \|\phi_H\|_{L^1}}   \cdot t\cdot \nu_H(A) \Big),\quad t\geq 0, A\in\mathscr{U}.
 \eeqnn
 Specially, we obtain the following {\it functional} central limit theorem for marked Hawkes processes, which extends the central limit theorem established in Karabash and Zhu \cite{KarabashZhu2015}.
 \begin{corollary}\label{Thm315}
 Assume that Condition~\ref{C1} holds and $\|\phi_H\|_{L^1}<1$, we have $\{\sqrt{T}\bar{N}_{H,t}^T (A):t\geq 0, A\in\mathscr{U}  \}$ converges weakly to $\{W_{H,t}(A) + W_{I,t}(A) + \sigma_Z B_Z(t) \cdot \nu_H(A):t\geq 0, A\in\mathscr{U} \}$ in the space $\mathbb{D}([0,\infty),\mathcal{S}(\mathbb{U}))$ as $T\to\infty$, where the terms in the limit are defined as before with $\lambda_I=\lambda $ and $\nu_I(du)=\nu_H(du)$. 
 Specially, denote by $N_\lambda(t):= N_\lambda([0,t],\mathbb{U})$ the marked Hawkes process, we have
 \beqnn
 \sqrt{T} \Big(\frac{1}{T}N_{\lambda}(Tt)-\frac{\lambda\cdot t}{1-\|\phi_H\|_{L^1}} \Big)\ar\to\ar \sqrt{\lambda}\cdot B_I(t) + \sqrt{\frac{\lambda\cdot \|\phi_H\|_{L^1}}{1-\|\phi_H\|_{L^1} }} \cdot B_H(t) + \sigma_Z B_Z(t) ,
 \eeqnn
 where $\{(B_I(t),B_H(t)):t\geq 0\} $ is a two-dimensional standard Brownian motion and 
 \beqnn
 \langle B_I , B_Z\rangle_t= \frac{\sqrt{\lambda}}{\sigma_Z} \frac{ 1}{1-\|\phi_H\|_{L^1} }
 \quad \mbox{and}\quad
  \langle B_H , B_Z\rangle_t= \frac{\sqrt{\lambda}}{\sigma_Z}  \frac{ \|\phi_H\|^{1/2}_{L^1}}{|1-\|\phi_H\|_{L^1}|^{3/2} }.
 \eeqnn
 	\end{corollary}

  \begin{example}{\rm (Multivariate Hawkes process with common intensity)
  	For the finite mark space $\mathbb{U}=\{1,\cdots,d \}$, the marked Hawkes point measure $N_\lambda(ds,du)$ reduces to a multivariate Hawkes process with common intensity. Indeed,  the $d$-dimensional point process $\{(N_1(t),\cdots,N_d(t)):t\geq 0\} $ defined by $N_i(t):= N_\lambda([0,t],\{i \})$ has intensity $\{w\cdot Z_\lambda(t):t\geq 0\}$, where  $w=(\nu_H(\{1\}),\cdots,\nu_H(\{d\}))$ and 
  	\beqnn
  	Z_\lambda(t)= \lambda+ \sum_{i=1}^d \int_0^t \phi(t-s,i) N_i(ds) .
  	\eeqnn
	The conditions in Corollary~\ref{Thm315} hold when   $\sup_{i=1}^d\|\phi(i)  \|_{\infty}<\infty$ and $ \sum_{i=1}^d \nu_H(\{i \}) \|\phi(i) \|_{L^1}<1$. In this case,  we have as $T\to \infty$,  
	\beqnn
	\sqrt{T}\Big( \frac{N_k(Tt)}{T}- \frac{\lambda \nu_H(\{k\})\cdot t }{1-  \sum_{i=1}^d\nu_H(\{i\})\|\phi(i)  \|_{L^1}} \Big) \to  c_k B_k(t) + c_m \nu_H(\{k\}) B_m(t) , \quad k\in\mathbb{U},
	\eeqnn
   weakly in the space $\mathbb{D}([0,\infty),\mathbb{R})$, where 
   	\beqnn
   |c_k|^2= \frac{\lambda \cdot \nu_H(\{k\})}{1- \sum_{i=1}^d\nu_H(\{i\})\|\phi(i) \|_{L^1}},
   \quad 
   |c_m|^2 =  \frac{\lambda \cdot \sum_{i=1}^d \nu_H(\{i\})\|\phi(i)  \|_{L^1}^2   }{|1- \sum_{i=1}^d \|\phi(i)  \|_{L^1}|^3},
   \eeqnn
   $\{(B_1(t),\cdots,B_d(t)): t\geq 0\}$ is a standard $d$-dimensional Brownian motion and $\{B_m(t):t\geq 0\}$ is the common Brownian motion satisfying that $ \langle B_k, B_m   \rangle_t =\frac{c_k}{c_m}\cdot \|\phi(k) \|_{L^1} \cdot t $ for $k\in\mathbb{U}$.

}
	
  \end{example}

   \section{Proofs of the auxiliary results}
\setcounter{equation}{0}
\medskip

 In this section, we give the proofs of Lemma~\ref{Thm303}, \ref{Thm401} and Proposition~\ref{Thm402}, \ref{Thm403}. The proofs are based on a new 
%
%
%
stochastic Volterra representation for the intensity process $\{Z(t):t\geq 0\}$. To this end, we first link marked Hawkes processes with immigration to Hawkes random measures as introduced in  \cite{HorstXu2019} through the following two-parameter processes: for any $t\geq 0$ and $u,u'\in\mathbb{U}$,
\beqlb\label{eqn5.01.01}
\mathcal{Z}(t,u') := Z(t)\mathbf{1}_{\{ u'\in\mathbb{U} \}},\quad
\mathcal{Z}_0(t,u') := \mu_0(t) \mathbf{1}_{\{ u'\in\mathbb{U} \}}\quad \mbox{and}\quad
\varPhi(t,u',u):= \phi(t,u)\mathbf{1}_{\{ u'\in\mathbb{U} \}}.
\eeqlb
From (\ref{eqn2.05}), it is easy to see that $ \{\mathcal{Z}(t,u') :t\geq 0,u'\in\mathbb{U} \}$ solves the following stochastic Volterra-Fredholm integral equation:
\beqlb\label{eqn5.01.02}
\mathcal{Z}(t,u')\ar=\ar  \mathcal{Z}_0(t,u') + \int_0^t \int_\mathbb{U} \varPhi(t-s,u',u ) N_I(ds,du) 
+\int_0^t \int_\mathbb{U}\int_0^{\mathcal{Z}(s-,u)} \varPhi(t-s, u',u)N_0(ds,du,dz).
\eeqlb
From (\ref{eqn2.03}) and (\ref{eqn5.01.02}), we see that  $N_H(dt,du')$  is a Hawkes random measure on $\mathbb{R}_+\times\mathbb{U}$ with $(\mathscr{F}_t)$-intensity $\{ \mathcal{Z}(t,u'):t\geq 0,u'\in\mathbb{U} \}$ and basis measure $\nu_H(du')$; see \cite[Definition~2.2]{HorstXu2019}. 
Applying Theorem~2.2 in \cite{Xu2018},  we see that $\{ \mathcal{Z}(t,u'):t\geq 0,u'\in\mathbb{U} \}$ also solves the following stochastic Volterra-Fredholm integral equation:
\beqlb\label{eqn5.01.03}
\mathcal{Z}(t,u')\ar=\ar  \int_0^t\int_\mathbb{U} \mathbf{R}(t-s,u',u) \mathcal{Z}_0(t,u) ds\nu_H(du) + \int_0^t \int_\mathbb{U}  \mathbf{R}(t-s,u',u)N_I(ds,du)\cr
\ar\ar +\mathcal{Z}_0(t,u')   +  \int_0^t \int_\mathbb{U}\int_0^{\mathcal{Z}(s-,u)} \mathbf{R}(t-s, u',u)\tilde{N}_0(ds,du,dz),
\eeqlb
where
\beqlb\label{eqn5.01.04}
\mathbf{R}(t, u',u)\ar=\ar \varPhi(t, u',u) + \int_0^t \int_\mathbb{U} \mathbf{R}(t-s, u',u'')\varPhi(s,u'',u)ds\nu_H(du'').
\eeqlb
 Integrating both sides of (\ref{eqn5.01.04}) with respect to the probability measure $\nu_H(du')$, we see that $\mathbf{R}(t,u',u)$ yields a decomposition of $R(t,u')$ introduced in (\ref{eqn4.05}) as
 \beqlb\label{eqn5.01.05}
 R(t,u)=\int_{\mathbb{U}} \mathbf{R}(t,u',u)\nu_H(du').
 \eeqlb
 
The following proposition yields the desired representation for the intensity process $\{Z(t):t\geq 0 \}$ in terms of the \textsl{martingale measure} $\tilde N_0$. The representation in terms of the martingale measure is key to our subsequent analysis. The proof follows  from integrating both sides of (\ref{eqn5.01.03}) with respect to the probability measure $\nu_H(du')$.
\begin{proposition}\label{Thm5.01.01}
	The intensity process $\{Z(t):t\geq 0 \}$ is the unique solution to the following stochastic Volterra integral equation:
	\beqlb\label{eqn5.01.07}
	Z(t)\ar=\ar \mu_0(t)+ \int_0^t R_H(t-s) \mu_0(t) ds + \int_0^t \int_{\mathbb{U}}  R(t-s,u)N_I(ds,du)\cr
	\ar\ar    +  \int_0^t \int_{\mathbb{U}} \int_0^{Z(s-)} R(t-s, u)\tilde{N}_0(ds,du,dz).
	\eeqlb
\end{proposition}

\subsection{Proofs of Lemma~\ref{Thm303}, \ref{Thm401} and Proposition~\ref{Thm402}, \ref{Thm403}}\label{ProofA}

 Armed with the representation (\ref{eqn5.01.07}), we can now give the proofs of  {Lemma~\ref{Thm303}, \ref{Thm401} and Proposition~\ref{Thm402}, \ref{Thm403}.}
 For any $T,t>0$, integrating both sides of (\ref{eqn5.01.07}) over the interval $(0,Tt]$ and changing the order of integration,   
 \beqlb\label{eqn5.02.01}
 \int_0^{Tt} Z(s)ds
 \ar=\ar  \int_0^{Tt}\mu_0(s)ds+  \int_0^{Tt} R_H(Tt-s) ds \int_0^s\mu_0(r) dr\cr
 \ar\ar    +  \int_0^t \int_{\mathbb{U}} \int_0^{Z(Ts-)} \Big(\int_0^{T(t-s)}R(r, u)dr  \Big)\tilde{N}_0(dTs,du,dz)\cr
 \ar\ar  + \int_0^t  \int_{\mathbb{U}}  \Big(\int_0^{T(t-s)}R(r, u)dr  \Big) N_I(dTs,du)
 \eeqlb
 and
 \beqlb\label{eqn5.02.02}
  \int_0^{Tt} Z(s)ds- \int_0^{Tt}\mathbf{E}[Z(s)]ds 
 \ar=\ar \int_0^{Tt} (\mu_0(s)-\mathbf{E}[ \mu_0(s)]) ds  +  \int_0^{Tt} R_H(Tt-s) ds \int_0^s(\mu_0(r)-\mathbf{E}[\mu_0(r)]) dr\cr
 \ar\ar  +   \int_0^t  \int_{\mathbb{U}} \int_0^{Z(Ts-)}\Big( \int_0^{T(t-s)}R(r,  u)dr \Big) \tilde{N}_0(dTs,du,dz)\cr
 \ar\ar  + \int_0^t \int_{\mathbb{U}} \Big(  \int_0^{T(t-s)}R(r,u)dr\Big) \tilde{N}_I	(dTs,du).
 \eeqlb
From Condition~\ref{C1}, we can see $\|\mu_0\|_{L^1} <\infty$ a.s.~so the first term on the right side of (\ref{eqn5.02.02}) is uniformly bounded.
Since $\|R_H\|_{L^1}<\infty$,  we also get the uniform boundedness of the  second term as well as the following lemma.

 \begin{lemma}\label{Thm502.01}
 Under Condition~\ref{C1}, for any $\kappa>0$ we have the following uniform convergence in probability: as $T\to\infty$
 \beqlb\label{eqn5.02.03}
 T^{-\kappa/2} \int_0^{Tt} (\mu_0(s)-\mathbf{E}[ \mu_0(s)]) ds  + T^{-\kappa/2}   \int_0^{Tt} R_H(Tt-s) ds \int_0^s(\mu_0(r)-\mathbf{E}[\mu_0(r)]) dr\overset{\mathbf{P}}\to 0.
 \eeqlb
 \end{lemma}

 Next, we consider the weak convergence of the last two stochastic Volterra integrals on the right side of (\ref{eqn5.02.02}). 
 From (\ref{eqn4.06}) and (\ref{eqn4.07}), for any $u\in\mathbb{U}$ we see that the integrand $\int_0^{T(t-s)}R(r,  u)dr $ increases to $\| R(u)\|_{L^1}$ as $T\to\infty$. 
 Thus, we may approximate these two stochastic Volterra integrals with the following two $(\mathscr{F}_{Tt})$-martingales, respectively:
 \beqlb
 M_{T,H}(t)\ar:=\ar  \int_0^t  \int_{\mathbb{U}} \int_0^{Z(Ts-)} \| R(u)\|_{L^1}\tilde{N}_0(dTs,du,dz),\label{eqn5.02.04} \\
 M_{T,I}(t)\ar:=\ar   \int_0^t  \int_{\mathbb{U}} \| R(u)\|_{L^1}\tilde{N}_I(dTs,du). \label{eqn5.02.05}
 \eeqlb
 We denote the error processes of the above approximations by $\{\varepsilon_{T,H}(t):t\geq 0\}$ and $\{\varepsilon_{T,I}(t):t\geq 0\}$, respectively. 
 They have the following representations:
 \beqlb
 \varepsilon_{T,H}(t)\ar:=\ar   \int_0^t  \int_{\mathbb{U}}  \int_0^{Z(Ts-)} \Big(\int_{T(t-s)}^\infty R(r, u)dr\Big) \tilde{N}_0(dTs,du,dz),\label{eqn5.02.06}\\
 \varepsilon_{T,I}(t)\ar:=\ar   \int_0^t  \int_{\mathbb{U}} \Big(\int_{T(t-s)}^\infty R(r,u)dr \Big) \tilde{N}_I(dTs,du).\label{eqn5.02.07}
 \eeqlb
 For any $\kappa>0$ and $t\geq 0$, let $ M^{\kappa}_{T,i}(t):=  T^{-\kappa/2}M_{T,i}(t)$ and  $\varepsilon^{\kappa}_{T,i}(t):=T^{-\kappa/2}\varepsilon_{T,i}(t)$ with $i\in\{H,I \}$. 
 Applying the  Burkholder-Davis-Gundy inequality for the martingale $\{M^\kappa_{T,H}(t):t\geq 0\} $, we get
 \beqlb\label{eqn5.02.08}
  \mathbf{E}\Big[\sup_{t\in[0,1]} | M^\kappa_{T,H}(t)|^2\Big] 
 \ar\leq\ar \frac{C}{T^{\kappa}}\mathbf{E}\Big[ \int_0^1 \int_{\mathbb{U}}  \int_0^{Z(Ts-)}
 \| R(u)\|_{L^1}^2 N_0(dTs,du,dz)  \Big]\cr
 \ar\leq\ar \frac{C}{T^{\kappa-1}} \int_0^1 \mathbf{E} [Z(Ts)]ds
 \int_{\mathbb{U}}  \| R(u)\|_{L^1}^2\nu_H(du).
 \eeqlb
 The integral in the last term above is uniformly bounded in $T$; see  (\ref{eqn4.07}) and Lemma~\ref{Thm201}. 
Similarly, we also can prove that $\mathbf{E}[\sup_{t\in[0,1]} | M^\kappa_{T,I}(t)|^2] \leq CT^{1-\kappa} $. 
These yield the following result.
 \begin{lemma}\label{Thm502.02}
 	Under Condition~\ref{C1},	for any $\kappa> 1$, both $\{M^\kappa_{T,H}(t):t\geq 0\}$ and $\{M^\kappa_{T,I}(t):t\geq 0\}$ converge weakly to $0$ in the space $\mathbb{D}([0,\infty),\mathbb{R})$ as $T\to\infty$.
 	\end{lemma}
 
The following lemma shows that $\{ \varepsilon^\kappa_{T,i}(t):t\geq 0 \}$ is weakly convergent to $0$ in the space $\mathbb{D}([0,\infty),\mathbb{R})$. 
The technical proof is given in Section~\ref{AppendixA}.
\begin{lemma}\label{Thm502.03}
	Under Condition~\ref{C1}, for any $\kappa\geq 1$,  both  $\{\varepsilon^\kappa_{T,H}(t):t\geq 0\}$ and $\{\varepsilon^\kappa_{T,I}(t):t\geq 0\}$ converge weakly to $0$ in the space $\mathbb{D}([0,\infty),\mathbb{R})$ as $T\to\infty$.
\end{lemma}
We are now ready to give the proofs of the auxiliary results.

\smallskip
\noindent\textit{Proof of Lemma~\ref{Thm303}.} Firstly, we have 
\beqlb
\frac{1}{T}\int_0^{Tt} Z(s)ds- \|R_I\|_{L^1}\cdot t \ar=\ar  \frac{1}{T}\Big(  \int_0^{Tt} Z(s)ds- \int_0^{Tt}\mathbf{E}[Z(s)]ds  \Big) \cr
\ar\ar +   \frac{1}{T}\int_0^{Tt}\mathbf{E}[Z(s)]ds - \lambda_I\cdot \|R_I\|_{L^1}\cdot t.
\eeqlb
The uniformly convergence of the second term on the right side of the equation above to $0$ follows from  (\ref{eqn3.14}). For the first term, applying Lemma~\ref{Thm502.01}, \ref{Thm502.02} and \ref{Thm502.03} with $\kappa=2$, we have 
\beqlb
\frac{1}{T}\Big(  \int_0^{Tt} Z(s)ds- \int_0^{Tt}\mathbf{E}[Z(s)]ds  \Big) \to 0, 
\eeqlb
weakly in the space $\mathbb{D}([0,\infty),\mathbb{R})$ and hence uniformly in probability on any bounded interval; see \cite[p.124]{Billingsley1999}.  This finishes the proof. 
\qed

\smallskip
\noindent\textit{Proof of Lemma~\ref{Thm401} and Proposition~\ref{Thm402}.}  
We give a detailed proof of Proposition~\ref{Thm402}. 
The proof of Lemma~\ref{Thm401} is similar but simpler. 
 For any $f \in B(\mathbb{U})$,  define 
 \beqlb
 W^f_T(t)\ar:=\ar \int_0^t\int_{\mathbb{U}}\int_0^{Z(Ts-)}\frac{f(u)}{\sqrt{T}} \tilde{N}_0(dTs,du,dz).
 \eeqlb
 In what follows we verify that $\{W^f_T(t):t\geq 0\}_{T\geq 0}$  satisfies the conditions of the Lindeberg-Feller theorem; see Theorem~3.22 in \cite[p.476]{JacodShiryaev2003}. 
 In this case, $\{W^f_T(t):t\geq 0 \}_{T\geq 0}$ converges weakly to a Brownian motion. 
 The condition~(3.23) in \cite[p.476]{JacodShiryaev2003} follows directly from the following statement:  for any $t\geq 0$, 
 \beqlb\label{eqn5.15.01}
 \mathbf{E}\Big[\sum_{s\leq t}|W^f_T(s)-W^f_T(s-)|^{2\alpha}\Big]\ar=\ar \mathbf{E}\Big[\int_0^t \int_{\mathbb{U}} \int_0^{Z(Ts-)}
 \frac{|f(u)|^{2\alpha}}{T^\alpha}
 N_0(dTs,du,dz)\Big]\cr
 \ar\leq\ar  \frac{C}{T^{\alpha-1}}\int_0^t \mathbf{E}[Z(Ts)]ds,
 \eeqlb
 which vanishes as $T\to\infty$; see Lemma~\ref{Thm201}.
 It remains to prove that condition [$\hat{\gamma}_5'$-$D$] in \cite[p.473]{JacodShiryaev2003} holds for $\{W^h_T(t):t\geq 0\}_{T\geq 0}$, i.e. as $T\to\infty$,
 \beqlb\label{eqn5.17.01}
 [W^f_T,W^f_T]_t\overset{\mathbf{P}}\to  \lambda_I\cdot  \| R_I\|_{L^1} \cdot \nu_H(|f|^2)\cdot  t.
 \eeqlb
 The quadratic variation of martingale $\{W^f_T(t):t\geq 0\}$ has the following  representation:
 \beqlb\label{eqn5.19.01}
 [W^f_T,W^f_T]_t\ar=\ar \int_0^{Tt} \int_{\mathbb{U}} \int_0^{Z(s-)}
 \frac{|f(u)|^2}{T} N_0(ds,du,dz).
 \eeqlb
 From (\ref{eqn3.14}),
 it is sufficient to prove that the following martingale
 \beqlb\label{eqn5.20.01}
 \mathcal{M}_{T}^f(t)\ar:=\ar [W^f_T,W^f_T]_t-\int_0^{Tt}Z(s)ds  \int_{\mathbb{U}} \frac{|f(u)|^2}{T}  \nu_H(du)\cr
 \ar=\ar \int_0^{Tt}  \int_{\mathbb{U}}  \int_0^{Z(s-)} \frac{|f(u)|^2}{T}  \tilde{N}_0(ds,du,dz)
 \eeqlb
  converges to $0$ weakly in the space $\mathbb{D}([0,\infty),\mathbb{R})$ as $T\to\infty$. The convergence of $\{ \mathcal{M}_{T}^f(t):t\geq 0\}_{T\geq0}$ to $0$ in the sense of finite-dimensional distributions can be seen as follows: from the Burkholder-Davis-Gundy inequality and Lemma~\ref{Thm201},
 \beqlb\label{eqn5.21.01}
 \mathbf{E}[|\mathcal{M}_{T}^f(t)|^{\alpha} ]\ar\leq\ar \frac{C}{T^\alpha} \mathbf{E}\Big[\Big| \int_0^{Tt}  \int_{\mathbb{U}}  \int_0^{Z(s-)}
 |f(u) |^4 N_0(ds,du,dz)\Big|^{\frac{\alpha}{2} }\Big]\cr
 \ar\leq\ar \frac{C}{T^\alpha} \mathbf{E}\Big[ \int_0^{Tt}  \int_{\mathbb{U}}  \int_0^{Z(s-)}
 |f(u) |^{2\alpha} N_0(ds,du,dz)\Big]\cr
 \ar\leq\ar \frac{C}{T^{\alpha-1}} \int_0^{t}  \mathbf{E}[Z(Ts)]ds ,
 \eeqlb
  which vanishes as $T\to\infty$. Now we show that the sequence $\{ \mathcal{M}_{T}^f (t):t\geq 0\}_{T\geq 0}$ is tight.
  To this end, we notice that the  sample paths  of $\{ \mathcal{M}_{T}^f (t):t\geq 0\}$ have total variation 
 \beqlb\label{eqn5.22.01}
 \mathrm{TV}(\mathcal{M}_{T}^f)(t)\ar=\ar \frac{1}{T}\int_0^{Tt} \int_{\mathbb{U}}  \int_0^{Z(s-)}
 |f(u)|^{2} N_0(ds,du,dz)+ \nu_H(|f|^2) \cdot \int_0^{t} Z(Ts)ds   ,
 \eeqlb
 which is finite almost surely. 
 From Theorem~3.36 in \cite[p.354]{JacodShiryaev2003}, it thus suffices to prove that the sequence of increasing process  $\{\mathrm{TV}(\mathcal{M}_{T,0})(t):t\geq 0\}_{T\geq 0}$ is $C$-tight. The $C$-tightness of the second integral on the right side of (\ref{eqn5.22.01}) follows directly from Kurtz's criterion or Aldous's criterion; see Theorem~6.8 in \cite{Walsh1986}. 
 For the first integral, from Theorem~3.37 in Jacod and Shiryeav \cite[p.354]{JacodShiryaev2003}, it suffices to prove that it converges to the   linearly increasing function $\lambda_I\cdot \|R_I\|_{L^1} \cdot \nu_H(|f|^2)\cdot t$
 in the sense of finite-dimensional distributions, which follows directly from (\ref{eqn5.19.01})-(\ref{eqn5.21.01}) and Theorem~\ref{Thm304}.
 Hence we have proved that there exists a Brawnian motion $\{W_{H,t}(f): t\geq 0\}$ with quadratic variation $\lambda_I\cdot \|R_I\|_{L^1} \cdot  \nu_H(|f|^2)\cdot t$ such that $\{W^f_T(t):t\geq 0\}$ converges to $ \{W_{H,t}(f):t\geq 0\}$ weakly in the space of $\mathbb{D}([0,\infty),\mathbb{R})$. 
 
  Using similar arguments, we can also prove that for any $f,f'\in B(\mathbb{U})$,
  \beqnn
  [ W_H(f),W_H(f') ]_t \overset{\mathbf{P}}= \lim_{T\to\infty}  [W^f_T,W^{f'}_T ]_t \ar=\ar \lim_{T\to\infty} \frac{1}{T}\int_0^{Tt} \int_{\mathbb{U}}  \int_0^{Z(s-)}
   f(u) f'(u)N_0(ds,du,dz) \cr
   \ar\ar\cr
   \ar\overset{\mathbf{P}}=\ar \lambda_I\cdot  \|R_I\|_{L^1} \int_{\mathbb{U}}f(u)f'(u) \nu_H(du)\times t. 
  \eeqnn
 This implies that $\{ W_{H,t}(A):t\geq 0, A\in \mathbb{U} \}$ is a continuous, worthy martingale measure on $\mathbb{U}$  with covariance measure $Q_H(ds,du,du')=ds\nu_H(du)\delta_u(du')$, where $\delta_u(du')$ is a Dirac measure on the point $u$.
 From  \cite[Proposition~2.10]{Walsh1986}, we can see that $\{ W_{H,t}(A):t\geq 0, A\in \mathbb{U} \}$ is a  Gaussian white noise on $\mathbb{U}$ with intensity $ \lambda_I\cdot \|R_I\|_{L^1}\cdot dt\nu_H(du)$.
 The independence of $\{ W_{H,t}(A):t\geq 0, A\in \mathbb{U} \}$ and $\{ W_{I,t}(A):t\geq 0, A\in \mathbb{U} \}$ follows directly from the fact that $N_0(ds,du,dz)$ and $N_1(ds,du)$ are independet. 
  \qed

 \smallskip
 \noindent\textit{Proof of Proposition~\ref{Thm403}.}  From Proposition~\ref{Thm302}, it suffices to prove that 
 \beqnn
 \frac{1}{\sqrt{T}}\Big(  \int_0^{Tt} Z(s)ds- \int_0^{Tt}\mathbf{E}[Z(s)]ds  \Big) \to \sum_{i\in\{H,I\}} \int_0^t \int_{\mathbb{U}}\|R(u)\|_{L^1} W_i(ds,du),
 \eeqnn
 weakly in the space $\mathbb{D}([0,\infty),\mathbb{R})$.
 From (\ref{eqn5.02.02}) and (\ref{eqn5.02.04})-(\ref{eqn5.02.07}), 
 \beqnn
 \frac{1}{\sqrt{T}}\Big(  \int_0^{Tt} Z(s)ds- \int_0^{Tt}\mathbf{E}[Z(s)]ds  \Big) \ar=\ar \frac{1}{\sqrt{T}}\int_0^{Tt} (\mu_0(s)-\mathbf{E}[ \mu_0(s)]) ds \cr
 \ar\ar  + \frac{1}{\sqrt{T}} \int_0^{Tt} R_H(Tt-s) ds \int_0^s(\mu_0(r)-\mathbf{E}[\mu_0(r)]) dr\cr
 \ar\ar\cr
 \ar\ar +  M^{1}_{T,H}(t) + M^{1}_{T,I}(t)+ \varepsilon^{1}_{T,H}(t)+ \varepsilon^{1}_{T,H}(t).
 \eeqnn
 Applying Lemma~\ref{Thm401} and Proposition~\ref{Thm402} to (\ref{eqn5.02.04})-(\ref{eqn5.02.05}), for any $ i\in\{ H, I\}$ we have $\{M^{1}_{T,i}(t):t\geq 0  \}$ converges to $\{ W_{i,t}(\|R(\cdot)\|_{L^1}):t\geq 0  \}$ weakly in $\mathbb{D}((0,\infty),\mathbb{R})$ as $T\to\infty$.
 From this and Lemma~\ref{Thm502.01}-\ref{Thm502.03}, we have as $T\to\infty$, 
 \beqnn
  \frac{1}{\sqrt{T}}\Big(  \int_0^{Tt} Z(s)ds- \int_0^{Tt}\mathbf{E}[Z(s)]ds  \Big)  \to W_{H,t}(\|R(\cdot)\|_{L^1})+W_{I,t}(\|R(\cdot)\|_{L^1})
 \eeqnn
 weakly in the space $\mathbb{D}([0,\infty),\mathbb{R})$.
  From the properties of Gaussian white noise,  it is easy to check that $\{ W_{H,t}(\|R(\cdot)\|_{L^1})+W_{I,t}(\|R(\cdot)\|_{L^1}):t\geq 0  \} $
 is a Brownian motion with quadratic variation as follows
 \beqnn
 \lambda_I\cdot \|R_I\|_{L^1} \int_{\mathbb{U}}\|R(u)\|_{L^1}^2 \nu_H(du) \cdot t + \lambda_I\cdot \int_{\mathbb{U}}\|R(u)\|_{L^1}^2\nu_I(du)\cdot t.
 \eeqnn
 The desired result (\ref{eqn4.08.01}) follows directly from this and (\ref{eqn4.07}). 
 \qed
 


 \subsection{Proof of Lemma~\ref{Thm502.03}}\label{AppendixA}

 In this section, we  prove the weak convergence of error processes $\{ \varepsilon^\kappa_{T,H} (t):t\geq 0\}_{T\geq 0}$ to $0$ for $\kappa=1$ under Condition~\ref{C1}.
 The proof for the sequence $\{ \varepsilon^\kappa_{T,I}(t):t\geq 0 \}_{T\geq 0}$ is similar but much simpler.
 
 For any $T>0$, let $\delta_T=T^{-\beta}$ for some $\beta\in(\frac{1}{\alpha}, 1-\frac{1}{2\theta_0})$.
 We split the error process  $\varepsilon^1_{T,H}(t)$ into the following two parts:
 \beqlb
 \varepsilon_{T,H,1}(t)
 \ar:=\ar \frac{1}{\sqrt{T}}  \int_0^{t-\delta_T} \int_{\mathbb{U}}  \int_0^{Z(Ts-)}
 \Big(\int_{T(t-s)}^\infty R(r,u)dr\Big)\tilde{N}_0(dTs,du,dz),  \label{eqnA.01} \\
 \varepsilon_{T,H,2}(t)
 \ar:=\ar \frac{1}{\sqrt{T}} \int_{t-\delta_T}^t \int_{\mathbb{U}}   \int_0^{Z(Ts-)}
\Big( \int_{T(t-s)}^\infty R(r,u)dr\Big)\tilde{N}_0(dTs,du,dz)).\label{eqnA.02}
 \eeqlb
 We first prove  that $\{\varepsilon_{T,H,1}(t):t\geq 0\}$ converges to $0$ uniformly in probability on any bounded interval as $T\to\infty$. The proof uses the following lemma.
 
 \begin{lemma}\label{ThmA01}
 	There exists a constant $C>0$ such that for any $t\geq 0$ and $i\in\{ H,I \}$,
 	\beqlb\label{eqnA.03}
 	\int_t^{\infty} R_i(s)ds  \leq  Ct^{-\theta_0}.
 	\eeqlb
 \end{lemma}
 \proof
 Let $\{X_k\}_{k\geq 1}$ be a sequence of i.i.d. random variables with probability density function $  \phi_H(t)/\|\phi_H\|_{L_1}$ and $N_G$  be a geometric random variable with parameter $1-   \|\phi_H\|_{L_1}>0$ independent of $\{X_k\}_{k\geq 1}$.
 From the one-to-one correspondence  between probability laws and their Laplace transforms, it is easy to see that the geometric summation $\sum_{k=1}^{N_G}X_k$ has the following probability density function
 \beqlb
 \frac{1-\|\phi_H\|_{L_1}}{\|\phi_H\|_{L_1}} \sum_{i=1}^\infty  \phi_H^{(*i)}(t)
 = \frac{1-\|\phi_H\|_{L_1}}{\|\phi_H\|_{L_1}} R_H(t) . \label{eqnA.03.1}
 \eeqlb
 From Condition~\ref{C1}, the geometric summation has finite $\theta_0$-th moment, i.e., $\mathbf{E}[|\sum_{k=1}^{N_G}X_k|^{\theta_0}]<\infty$.
 From Markov's inequality, there exits a constant $C>0$ such that for  any $t>0$,
 \beqnn
 \mathbf{P}\Big\{ \sum_{k=1}^{N_G}X_k  >t\Big\}
 \leq \mathbf{E}\Big[\Big|\sum_{k=1}^{N_G}X_k\Big|^{\theta_0}\Big]t^{-\theta_0}
 \leq Ct^{-\theta_0}.
 \eeqnn
 Taking this estimate back into (\ref{eqnA.03.1}), we have
 \beqnn
 \int_t^\infty  R_H(s)ds \leq C\mathbf{P}\Big\{ \sum_{k=1}^{N_G}X_k  >t\Big\}\leq Ct^{-\theta_0}.
 \eeqnn
 This yields the desired result for $i=H$.
 For the case $i=I$, from (\ref{eqn2.10}) we have
 \beqnn
 R_I(t)\ar=\ar \phi_I(t)+\frac{\|\phi_H\|_{L^1}\|\phi_I\|_{L^1}} {1-\|\phi_H\|_{L^1}} \int_0^t \frac{1-\|\phi_H\|_{L^1}}{\|\phi_H\|_{L^1}} R_H(t-s) \times \frac{\phi_I(s)}{\|\phi_I\|_{L^1}}ds.
 \eeqnn
 Let $Y$ be an $\mathbb{R}_+$-valued random variable with  probability density $\phi_I(t)/\|\phi_I\|_{L^1}$ and independent of  $\sum_{k=1}^{N_G}X_k$. Thus, the convolution in the above equation equals the density function of $\sum_{k=1}^{N_G}X_k+ Y$ whose $\theta_0$-th moment is finite. As before, we also have
 \beqnn
 \int_t^\infty R_I(s)ds \leq  t^{\theta_0}\int_t^\infty |\frac{s}{t}|^{\theta_0}\phi_I(s)ds + C \mathbf{P}\Big\{ \sum_{k=1}^{N_G}X_k+ Y \geq t\Big\} \leq C t^{-\theta_0}.
 \eeqnn 
 \qed

 \begin{proposition}\label{ThmA02}
 	We have $  \sup_{t\in[0,1]} |\varepsilon_{T,H,1}(t)|\to 0$ in probability as $T\to\infty$
 \end{proposition}
 \proof From the definition of $\beta$ and $\delta_T$, we see that $T(t-s)\geq T^{1-\beta}$ for any $s\in[0,t-\delta_T]$ and from (\ref{eqnA.01}),  
 \beqnn
 |\varepsilon_{T,H,1}(t)|
 \ar\leq \ar  \sqrt{T} \int_0^{t-\delta_T} Z(Ts)ds \int_{\mathbb{U}}
 \int_{T^{1-\beta}}^\infty R(r,u)dr \nu_H(du)  \cr
 \ar\ar  + \frac{1}{\sqrt{T}} \int_0^{t-\delta_T} \int_{\mathbb{U}}  \int_0^{Z(Ts-)}
\Big( \int_{T^{1-\beta}}^\infty R(r,u)dr \Big)N_0(dTs,du,dz).
 \eeqnn
 Taking expectation on the both sides of this inequality, it follows from Lemma~\ref{Thm201} that
 \beqnn
 \mathbf{E}\Big[\sup_{t\in[0,1]}|\varepsilon_{T,H,1}(t)|\Big]
 \ar\leq\ar 2\sqrt{T}\int_0^{1} \mathbf{E}[Z(Ts)]ds
 \int_{T^{1-\beta}}^\infty R_H(r)dr
 \leq C\sqrt{T}\int_{T^{1-\beta}}^\infty R_H(r)dr.
 \eeqnn
 From Lemma~\ref{ThmA01} and $\beta<1-\frac{1}{2\theta_0}$, we have as $T\to\infty$,
 \beqnn
 \mathbf{E}\Big[\sup_{t\in[0,1]}|\varepsilon_{T,H,1}(t)|\Big]\leq C T^{1/2- (1-\beta)\theta_0}\to 0.
 \eeqnn 
 \qed
 
 Now we prove that the sequence $\{ \varepsilon_{T,H,2}(t):t\geq 0 \}_{T\geq 1}$ converges weakly to $0$ as $T\to\infty$. It suffices to prove that it is tight in the space $\mathbb{D}([0,\infty),\mathbb{R})$ and converges $0$ in the sense of finite-dimensional distributions.
 
 \begin{proposition}\label{ThmA03}
 	The process $\{ \varepsilon_{T,H,2}(t):t\geq 0 \}$ converges to $0$  in the sense of finite-dimensional distributions as $T\to \infty$, i.e.,  for any $t\geq 0$,
 	\beqnn
 	\mathbf{E}\big[ |\varepsilon_{T,H,2}(t)|^{2}\big]\to 0.
 	\eeqnn
 \end{proposition}
 \proof From the  Burkholder-Davis-Gundy inequality, Lemma~\ref{Thm201} and (\ref{eqn4.07}), we have 
 \beqnn
 \mathbf{E}[ |\varepsilon_{T,H,2}(t)|^{2}]\ar\leq\ar \frac{C}{T}\mathbf{E}\Big[  \Big|   \int_{t-\delta_T}^t  \int_{\mathbb{U}}  \int_0^{Z(Ts-)}
 \Big|\int_{0}^\infty R(r,u)dr\Big|^2 N_0(dTs,du,dz)  \Big]\cr
 \ar\leq\ar C  \int_{t-\delta_T}^t   \mathbf{E}[Z(Ts)]ds  \int_{\mathbb{U}}
 \Big|\int_{0}^\infty R(r,u)dr\Big|^2 \nu_H(du) \leq CT^{-\beta},
 \eeqnn
 which vanishes as $T\to \infty$.
 \qed
 
 For the tightness of  $\{ \varepsilon_{T,H,2}(t):t\geq 0 \}_{T\geq 0}$ , by Theorem~8.8 in \cite[p.139]{EthierKurtz1986},
 it suffices to  prove that there exist two constants $\gamma>1$ and $C>0$ such that for any $T,h>0$,
 \beqlb\label{eqnA.04}
 \sup_{t\in[0,1]}\mathbf{E}\Big[ |\triangle_h\varepsilon_{T,H,2}(t)|^\alpha\times |\triangle_h\varepsilon_{T,H,2}(t+h)|^\alpha \Big]\leq Ch^{\gamma},
 \eeqlb
 where $\triangle_h f(t):=f(t)-f(t-h)$.
 As a preparation, we first give some high-order moment estimates for the intensity process $\{Z(t):t\geq 0 \}$ and the stochastic integral driven by $\tilde{N}_0(ds,du,dz)$.
 \begin{lemma}\label{ThmA04}
 	Under Condition~\ref{C1}, there exits a constant $C>0$ such that 
 	\begin{enumerate}
 		\item[(1)]  $\mathbf{E}[|Z(t)|^{2\alpha}]\leq C$ for any $t\geq 0$;
 		
 		\item[(2)] For any $T>0$, $0\leq t'\leq t''$ and any measurable function $f(t_1,t_2,t_3,u)$ on $\mathbb{R}_+^3\times \mathbb{U}$ satisfying that $\nu_0(|f(t_1,t_2,t_3,\cdot)|^{2\alpha})$ is locally bounded on $\mathbb{R}_+^3$, 
 		\beqlb\label{eqnA.05.01}
 		\lefteqn{\mathbf{E}\Big[\Big| \int_{t'}^{t''} \int_{\mathbb{U}}\int_0^{Z(Ts-)} f(t',t'',s,u)\tilde{N}_0(dTs,du,dz)\Big|^{2\alpha} \Big] } \ar\ar\cr
 		\ar\leq\ar  CT\int_{t'}^{t''}  ds \int_\mathbb{U} |f(t',t'',s,u)|^{2\alpha}\nu_H( dy)\cr
 		\ar\ar + CT^\alpha \int_{\mathbb{U}}   \nu_H(du) \int_{t'}^{t''}   |f(t',t'',r,u)|dr\cdot     \Big| \int_{t'}^{t''}  |f(t',t'',s,u)|^{\frac{2\alpha-1}{\alpha-1}}ds\Big|^{\alpha-1} . 
 		\eeqlb
 		\end{enumerate} 
 \end{lemma}
 \proof We first prove that $\sup_{t\geq 0}\mathbf{E}[|Z(t)|^2]<\infty$. 
 From (\ref{eqn5.01.07}) and the Cauchy-Schwarz inequality,
 \beqlb\label{eqnA.05.2}
 \mathbf{E}[|Z(t)|^2]
 \ar\leq\ar C\mathbf{E}[|\mu_0(t)|^2]+ C\mathbf{E}\Big[\Big|\int_0^t R_H(t-s)\mu_0(s)ds\Big|^2\Big]\cr
 \ar\ar +C\mathbf{E}\Big[\Big|\int_0^t  \int_{\mathbb{U}}  R(t-s,u)\tilde{N}_I(ds,du)\Big|^2\Big] \cr
 \ar\ar + C\mathbf{E}\Big[ \Big|\int_0^t \int_{\mathbb{U}} \int_0^{Z(s-)} R(t-s,u)\tilde N_0(ds,du,dz)\Big|^2\Big].
 \eeqlb
 The uniform boundedness of the first term on the right side of the above inequality follows from Condition~\ref{C1}. For the second term, from H\"older's inequality and (\ref{eqn3.06}), 
 \beqlb\label{eqnA.05.03}
 \mathbf{E}\Big[\Big|\int_0^t R_0(t-s)\mu_0(s)ds\Big|^2\Big]\ar\leq\ar \int_0^t R_H(r)dr \cdot \int_0^t R_H(t-s)\mathbf{E}[|\mu_0(s)|^2]ds\leq C \|R_H\|^2_{L^1}<\infty.
 \eeqlb
 Next, we prove that the last expectation  in (\ref{eqnA.05.2}) is uniformly bounded; the second-to-last term can be handled in the same way.
 From the Burkholder-Davis-Gundy inequality, Lemma~\ref{Thm201}, the Cauchy-Schwarz inequality and (\ref{eqn4.06})-(\ref{eqn4.07}),
 \beqlb \label{eqnA.05.04}
 \mathbf{E}\Big[ \Big|\int_0^t \int_{\mathbb{U}} \int_0^{Z(s-)} R(t-s,u)\tilde N_0(ds,du,dz)\Big|^{2}\Big] 
 \ar\leq\ar C \mathbf{E}\Big[\int_0^t \int_{\mathbb{U}} \int_0^{Z(s-)} |R(t-s,u)|^2 N_0(ds,du,dz)\Big]\cr  
 \ar\leq\ar  C  \int_0^t\mathbf{E}[Z(s)]ds  \int_{\mathbb{U}}  |R(t-s,u)|^2 \nu_H(du)\cr
 \ar\leq\ar C \int_{\mathbb{U}}  \|\phi(u)\|_\infty \cdot \|\phi(u)\|_{L^1} \nu_H(du)\cr
 \ar\leq \ar C \int_{\mathbb{U}}  [ \|\phi(u)\|_\infty^2 +  \|\phi(u)\|_{L^1}^2] \nu_H(du)<\infty.
 \eeqlb
 Taking these estimates back into (\ref{eqnA.05.2}), we can get  $\sup_{t\geq 0}\mathbf{E}[|Z(t)|^2]<\infty$. 
 Let us now prove (\ref{eqnA.05.01}). 
 From the Burkholder-Davis-Gundy inequality and the Cauchy-Schwarz inequality, 
 \beqlb\label{eqnA.07.00}
	\lefteqn{\mathbf{E}\Big[\Big| \int_{t'}^{t''} \int_{\mathbb{U}}\int_0^{Z(Ts-)} f(t',t'',s,u)\tilde{N}_0(dTs,du,dz)\Big|^{2\alpha} \Big] } \ar\ar\cr
\ar\leq\ar\mathbf{E}\Big[\Big| \int_{t'}^{t''} \int_{\mathbb{U}}\int_0^{Z(Ts-)} |f(t',t'',s,u)|^2N_0(dTs,du,dz)\Big|^{\alpha} \Big] \cr
\ar\leq\ar CT^\alpha \mathbf{E}\Big[\Big| \int_{t'}^{t''} Z(Ts)ds\int_{\mathbb{U}} |f(t',t'',s,u)|^2\nu_H(du)\Big|^{\alpha} \Big]\cr
\ar\ar +C \mathbf{E}\Big[\Big| \int_{t'}^{t''} \int_{\mathbb{U}}\int_0^{Z(Ts-)} |f(t',t'',s,u)|^2\tilde{N}_0(dTs,du,dz)\Big|^{\alpha} \Big].
 \eeqlb
 Applying Jensen's inequality and H\"older's inequality to the first term on the right side of the last inequality above, we can see that it can be bounded by
 \beqlb\label{eqnA.07.01}
  \lefteqn{CT^\alpha\int_{\mathbb{U}} \mathbf{E}\Big[\Big| \int_{t'}^{t''} Z(Ts) |f(t',t'',s,u)|^2ds\Big|^{\alpha} \Big]\nu_H(du)}\ar\ar\cr
  \ar\leq\ar CT^\alpha\int_{\mathbb{U}} \nu_H(du) \int_{t'}^{t''}\mathbf{E}[ |Z(Ts)|^\alpha ] |f(t',t'',r,u)|dr \cdot \Big| \int_{t'}^{t''}   |f(t',t'',s,u)|^\frac{2\alpha-1}{\alpha-1}ds\Big|^{\alpha-1}\cr
  \ar\leq\ar CT^\alpha\int_{\mathbb{U}}\nu_H(du) \int_{t'}^{t''}  |f(t',t'',r,u)|dr \cdot   \Big| \int_{t'}^{t''}   |f(t',t'',s,u)|^\frac{2\alpha-1}{\alpha-1}ds\Big|^{\alpha-1} .
 \eeqlb
 Applying the Burkholder-Davis-Gundy inequality again to the last term in (\ref{eqnA.07.00}), we have 
 \beqnn
\lefteqn{ \mathbf{E}\Big[\Big| \int_{t'}^{t''} \int_{\mathbb{U}}\int_0^{Z(Ts-)} |f(t',t'',s,u)|^2\tilde{N}_0(dTs,du,dz)\Big|^{\alpha} \Big]}\ar\ar\cr
\ar\leq\ar C\mathbf{E}\Big[\Big| \int_{t'}^{t''} \int_{\mathbb{U}}\int_0^{Z(Ts-)} |f(t',t'',s,u)|^4N_0(dTs,du,dz)\Big|^{\alpha/2} \Big]\cr
 \ar\leq\ar C\mathbf{E}\Big[ \int_{t'}^{t''} \int_{\mathbb{U}}\int_0^{Z(Ts-)} |f(t',t'',s,u)|^{2\alpha}N_0(dTs,du,dz) \Big]\cr
  \ar\leq\ar CT\int_{t'}^{t''} \mathbf{E}[ Z(Ts)]ds \int_{\mathbb{U}}|f(t',t'',s,u)|^{2\alpha}\nu_H(du) \cr
  \ar\leq\ar  CT\int_{t'}^{t''} ds \int_{\mathbb{U}}|f(t',t'',s,u)|^{2\alpha}\nu_H(du).
 \eeqnn
 Taking this and (\ref{eqnA.07.01}) back into (\ref{eqnA.07.00}), we get the desired result (\ref{eqnA.05.01}).  
 
 We now prove the first result. Using similar arguments as the ones leading to (\ref{eqnA.05.2}) and (\ref{eqnA.05.03}), we also get 
 \beqnn
 \mathbf{E}[|Z(t)|^{2\alpha}]
 \ar\leq\ar C + C\mathbf{E}\Big[ \Big|\int_0^t \int_{\mathbb{U}} \int_0^{Z(s-)} R(t-s,u)\tilde N_0(ds,du,dz)\Big|^{2\alpha}\Big].
 \eeqnn
 The uniform boundedness of the last term in the above inequality follows from  (\ref{eqnA.05.01}) with $T=1$, $t'=0$, $t''=t$ and $f(0,t,s,u)= R(t-s,u)$. Indeed,  from (\ref{eqn4.06})-(\ref{eqn4.07}) and the Cauchy-Schwarz inequality,  
 \beqnn
 \lefteqn{ \mathbf{E}\Big[ \Big|\int_0^t \int_{\mathbb{U}} \int_0^{Z(s-)} R(t-s,u)\tilde N_0(ds,du,dz)\Big|^{2\alpha}\Big]}\ar\ar\cr
 \ar\leq\ar  C\int_0^t  ds \int_\mathbb{U} |R(s,u)|^{2\alpha}\nu_H(du)
 + C \int_{\mathbb{U}}   \nu_H(du) \int_0^t R(r,u)dr     \Big| \int_0^t |R(s,u)|^{\frac{2\alpha-1}{\alpha-1}}ds\Big|^{\alpha-1} \cr
 \ar\leq\ar C \int_\mathbb{U} \|R(u)\|_{\infty}^{\alpha}\|R(u)\|_{L^1}^{\alpha}\nu_H(du)
 \leq C \int_\mathbb{U} \left[ \|R(u)\|_{L^1}^{2\alpha}+\| R(u) \|_\infty^{2\alpha}\right]  \nu_H(du)  <\infty .
 \eeqnn
 \qed
 
 \begin{proposition}\label{ThmA05}
 	There exists a constant $C>0$ such that for any $T\geq 0$ and $\kappa\in [1, \alpha]$,
 	\beqlb\label{eqnA.04.1}
 	\sup_{t\geq 0}\mathbf{E}\Big[ |\varepsilon_{T,H,2}(t)|^{2\kappa}\Big] \leq C|\delta_T|^\kappa.
 	\eeqlb
 	Moreover, the inequality (\ref{eqnA.04}) holds with $\gamma=\alpha$ for any $h\geq \delta_T/2$.
 \end{proposition}
 \proof  It is easy to see that the second result follows directly from the first one, i.e., by the Cauchy-Schwarz inequality, for any $h\geq \delta_T/2$,
 \beqnn
 \mathbf{E}\Big[ |\triangle_h\varepsilon_{T,H,2}(t)|^\alpha\times |\triangle_h\varepsilon_{T,H,2}(t+h)|^\alpha \Big]
 \ar\leq\ar C\mathbf{E}\Big[ |\varepsilon_{T,H,2}(t)|^{2\alpha}\Big] + C\mathbf{E}\Big[ |\varepsilon_{T,H,2}(t-h)|^{2\alpha} \Big]\cr
 \ar\ar\cr
 \ar\ar +C \mathbf{E}\Big[ |\varepsilon_{T,H,2}(t+h)|^{2\alpha}\Big]\leq C|\delta_T|^\alpha \leq Ch^\alpha.
 \eeqnn
 We prove the ineqaulity (\ref{eqnA.04.1}) for $\kappa=\alpha$; the general case can be proved in the same way.
 Applying Lemma~\ref{ThmA04}(2) with $t'=t-\delta_T$, $t''=t$ and $f(t',t'',s,u)= \int_{T(t-s)}^\infty R(r,u)dr$, we have 
 \beqnn
  \lefteqn{\mathbf{E}[ |\varepsilon_{T,H,2}(t)|^{2\alpha}]\leq\frac{C}{T^{\alpha-1}}  \int_{t-\delta_T}^t ds \int_{\mathbb{U}} \Big| \int_{T(t-s)}^\infty R(r,u)dr\Big|^{2\alpha}\nu_H(du)}\ar\ar\cr
  \ar\ar  + C\int_{\mathbb{U}}   \nu_H(du) \int_{t-\delta_T}^t \Big(  \int_{T(t-r)}^\infty R(z,u)dz\Big)  dr  \cdot   \Big| \int_{t-\delta_T}^t\Big| \int_{T(t-s)}^\infty R(r,u)dr\Big|^{\frac{2\alpha-1}{\alpha-1}}ds\Big|^{\alpha-1}.
 \eeqnn
 From (\ref{eqn4.07}), (\ref{eqn3.07}) and $\beta<1$, the first term on the right side of this inequality can be bounded by
 \beqnn
  \frac{C}{T^{\alpha-1}}  \int_{t-\delta_T}^t \int_\mathbb{U}  \|R(u)\|_{L^1}^{2\alpha} ds\nu_H(du)\leq  \frac{C \delta_T}{T^{\alpha-1}}= C|\delta_T|^{1+\frac{\alpha -1}{\beta}}\leq C|\delta_T|^\alpha.
 \eeqnn
 Similarly, the second term can be bounded by
 \beqnn
 C\int_{\mathbb{U}}   \nu_H(du) \int_{t-\delta_T}^t  \|R(u)\|_{L^1} ds     \Big| \int_{t-\delta_T}^t\|R(u)\|_{L^1} ^{\frac{2\alpha-1}{\alpha-1}}ds\Big|^{\alpha-1} \leq C |\delta_T|^\alpha \int_{\mathbb{U}}   \|R(u)\|_{L^1}^{{2\alpha} }  \nu_H(du) \leq  C| \delta_T|^\alpha  .
 \eeqnn
 Altogether, we obtain the desired result.
 \qed

 Now we are going to prove that (\ref{eqnA.04}) also holds for $h< \delta_T/2$. In this case, $t-h-\delta_T<t-\delta_T<t+h-\delta_T<t-h< t<t+h $, which suggests to decompose $\varepsilon_{T,H,2}(t+h)$ as:
 \beqnn
 \varepsilon_{T,H,2}(t+h)\ar=\ar  \frac{1}{\sqrt{T}} \int_{t}^{t+h} \int_{\mathbb{U}}   \int_0^{Z(Ts-)}
\Big( \int_{T(t+h-s)}^\infty R(r,u)dr\Big) \tilde{N}_0(dTs,du,dz)\cr
 \ar\ar +  \frac{1}{\sqrt{T}} \int_{t-h}^{t} \int_{\mathbb{U}}  \int_0^{Z(Ts-)}
 \Big( \int_{T(t+h-s)}^\infty R(r,u)dr\Big)\tilde{N}_0(dTs,du,dz)\cr
 \ar\ar +  \frac{1}{\sqrt{T}} \int_{t+h-\delta_T}^{t-h} \int_{\mathbb{U}}   \int_0^{Z(Ts-)}
\Big( \int_{T(t+h-s)}^\infty R(r,u)dr\Big) \tilde{N}_0(dTs,du,dz).
 \eeqnn
 We can decompose $\varepsilon_{T,H,2}(t)$ and $\varepsilon_{T,H,2}(t-h)$ in the same way. Thus,
 \beqnn
 \triangle_h\varepsilon_{T,H,2}(t+h)
 \ar=\ar \sum_{i=1}^4 I_{T,i}(t,h)\quad \mbox{and}\quad \triangle_h\varepsilon_{T,H,2}(t)=\sum_{j=5}^8 I_{T,j}(t,h),
 \eeqnn
 where
 \beqnn
 I_{T,1}(t,h)\ar=\ar  \frac{1}{\sqrt{T}} \int_{t}^{t+h} \int_{\mathbb{U}}  \int_0^{Z(Ts-)}
 \Big(\int_{T(t+h-s)}^\infty R(r,u)dr\Big) \tilde{N}_0(dTs,du,dz),\cr
 I_{T,2}(t,h)\ar=\ar -  \frac{1}{\sqrt{T}} \int_{t-h}^{t} \int_{\mathbb{U}}  \int_0^{Z(Ts-)}
 \Big(\int_{T(t-s)}^{T(t+h-s)} R(r,u)dr\Big) \tilde{N}_0(dTs,du,dz),\cr
 I_{T,3}(t,h)\ar=\ar -  \frac{1}{\sqrt{T}} \int_{t+h-\delta_T}^{t-h} \int_{\mathbb{U}}  \int_0^{Z(Ts-)}
 \Big(\int_{T(t-s)}^{T(t+h-s)} R(r,u)dr\Big) \tilde{N}_0(dTs,du,dz),\cr
 I_{T,4}(t,h)\ar=\ar - \frac{1}{\sqrt{T}} \int_{t-\delta_T}^{t+h-\delta_T} \int_{\mathbb{U}}  \int_0^{Z(Ts-)}
 \Big(\int_{T(t-s)}^\infty R(r,u)dr\Big)\tilde{N}_0(dTs,du,dz)
 \eeqnn
 and
 \beqnn
 I_{T,5}(t,h)\ar=\ar  \frac{1}{\sqrt{T}} \int_{t-h}^t \int_{\mathbb{U}}  \int_0^{Z(Ts-)}
 \Big( \int_{T(t-s)}^\infty R(r,u)dr\Big) \tilde{N}_0(dTs,du,dz),\cr
 I_{T,6}(t,h)\ar=\ar  -  \frac{1}{\sqrt{T}} \int_{t+h-\delta_T}^{t-h} \int_{\mathbb{U}}  \int_0^{Z(Ts-)}
 \Big(  \int_{T(t-h-s)}^{T(t-s)} R(r,u)dr\Big) \tilde{N}_0(dTs,du,dz),\cr
 I_{T,7}(t,h)\ar=\ar  - \frac{1}{\sqrt{T}} \int_{t-\delta_T}^{t+h-\delta_T} \int_{\mathbb{U}}  \int_0^{Z(Ts-)}
 \Big(  \int_{T(t-h-s)}^{T(t-s)} R(r,u)dr\Big) \tilde{N}_0(dTs,du,dz),\cr
 I_{T,8}(t,h)\ar=\ar  -  \frac{1}{\sqrt{T}} \int_{t-h-\delta_T}^{t-\delta_T} \int_{\mathbb{U}}  \int_0^{Z(Ts-)}
 \Big(  \int_{T(t-h-s)}^\infty R(r,u)dr\Big) \tilde{N}_0(dTs,du,dz).\cr
 \eeqnn
 From the Cauchy-Schwarz inequality,
 \beqnn
 \mathbf{E}\Big[ |\triangle_h\varepsilon_{T,H,2}(t)|^\alpha\times |\triangle_h\varepsilon_{T,H,2}(t+h)|^\alpha \Big]
 \ar\leq\ar C\sum_{i=1}^4\sum_{j=5}^8 \mathbf{E}\big[|I_{T,i}(T,h)|^\alpha \cdot |I_{T,j}(T,h)|^\alpha\big].
 \eeqnn
 Thus, it suffices to prove that each expectation in the sum can be bounded by $Ch^{\gamma}$ for some constants $C>0$ and $\gamma>1$ independent of $T$, $t$ and $h$. The following three results establish such bounds and hence finish the proof of Lemma~\ref{Thm502.03}. 
 
 \begin{proposition}\label{ThmA06}
 	There exists a constant $C>0$  such that  for any $t\in[0,1]$ any $h>0$, and all $i\in\{1,2,4\}$ and $j\in\{5,7,8\}$,
 	\beqnn
	\mathbf{E}\big[|I_{T,i}(t,h)|^\alpha \cdot |I_{T,j}(t,h)|^\alpha\big]\leq Ch^{\alpha\wedge \frac{3}{2}}.
 	\eeqnn
 \end{proposition}
 \proof We just consider the case $i=1$ and $j=5$. All other cases can be proved similarly.
 Using the similar arguments as in the proof of Proposition~\ref{ThmA05}, there exists a constant $C>0$ such that
 \beqlb\label{eqnA.05}
 \mathbf{E}[|I_{T,1}(t,h)|^{2\alpha}]+ \mathbf{E}[|I_{T,5}(t,h)|^{2\alpha}]\leq C\frac{(Th)^{\alpha}+Th}{T^\alpha}.
 \eeqlb
 When $Th\geq 1$, by H\"older's inequality,
 \beqnn
 \mathbf{E}\big[|I_{T,1}(t,h)|^\alpha \cdot |I_{T,5}(t,h)|^\alpha\big]\ar\leq\ar  \Big(\mathbf{E}\big[|I_{T,1}(t,h)|^{2\alpha}\big] \cdot  \mathbf{E}\big[|I_{T,5}(t,h)|^{2\alpha}\big]\Big)^{1/2} \leq Ch^\alpha.
 \eeqnn
 Now we consider the case $Th< 1$.
 For any $\xi\geq 0$, define
 \beqnn
 J_1(\xi)\ar:=\ar \int_{t}^\xi \int_{\mathbb{U}}  \int_0^{Z(Ts-)}\mathbf{1}_{\{s\in(t,t+h]\}}
 \Big(\int_{T(t+h-s)}^\infty R(r,u)dr\Big) \tilde{N}_0(dTs,du,dz),\cr
 J_2(\xi)\ar:=\ar \int_{t-h}^\xi \int_{\mathbb{U}}  \int_0^{Z(Ts-)} \mathbf{1}_{\{s\in(t-h,t]\}}
 \Big( \int_{T(t-s)}^\infty R(r,u)dr\Big) \tilde{N}_0(dTs,du,dz),
 \eeqnn
 which are two $(\mathscr{F}_{T\cdot})$-martingales for any $t$ and $h$ fixed. 
 Using the tower property of conditional expectation conditioning on $\mathcal{F}_{Tt}$, we have
 \beqlb\label{eqnA.06}
 \mathbf{E}\big[|I_{T,1}(t,h)|^\alpha \cdot |I_{T,5}(t,h)|^\alpha\big]
 \ar=\ar\frac{1}{T^\alpha}\mathbf{E}\big[|J_1(t+h)|^\alpha|J_2(t)|^\alpha\big]\cr
 \ar=\ar\frac{1}{T^\alpha} \mathbf{E}\big[|J_2(t)|^\alpha \cdot \mathbf{E}_{\mathscr{F}_{Tt}}[|J_1(t+h)|^\alpha]\big].
 \eeqlb
 As  before, from the Burkholder-Davis-Gundy inequality and (\ref{eqn4.07}), there exists a constant $C>0$ such that
 \beqnn
 \mathbf{E}_{\mathscr{F}_{Tt}}[|J_1(t+h)|^\alpha]
 \ar\leq\ar CT\int_{t}^{t+h} \mathbf{E}_{\mathscr{F}_{Tt}}[Z(Ts)]ds
 \eeqnn
 and
 \beqnn
 \mathbf{E}[|J_1(t+h)|^\alpha|J_2(t)|^\alpha]\ar\leq \ar CT\int_{t}^{t+h} \mathbf{E}\big[|J_2(t)|^\alpha \cdot \mathbf{E}_{\mathscr{F}_{Tt}}[Z(Ts)]\big]ds\cr
 \ar\leq\ar CT\int_{t}^{t+h}
 \mathbf{E}[|J_2(t)|^\alpha Z(Ts-)]ds.
 \eeqnn
 By H\"older's inequality, Lemma~\ref{ThmA04}(2) and (\ref{eqnA.05}),
 \beqnn
 \mathbf{E}[|J_1(t+h)|^\alpha|J_2(t)|^\alpha]\ar\leq\ar CT\int_{t}^{t+h}
 \sqrt{\mathbf{E}[|J_2(t)|^{2\alpha}]\cdot\mathbf{E}[ |Z(Ts)|^2]}ds\cr
 \ar\leq\ar  CTh
 \sqrt{\mathbf{E}[|J_2(t)|^{2\alpha}]}\leq CT^{3/2}h^{3/2}.
 \eeqnn
 Taking this back into (\ref{eqnA.06}), we have
 \beqnn
 \mathbf{E}\big[|I_{T,1}(t,h)|^\alpha \cdot |I_{T,5}(t,h)|^\alpha\big]\leq CT^{3/2-\alpha}h^{3/2}\leq \left\{
 \begin{array}{ll}
 	Ch^{3/2}, & \mbox{if }\alpha\geq 3/2;\vspace{7pt}\\
 	C h^\alpha, & \mbox{if }1<\alpha<3/2.
 \end{array}
 \right.
 \eeqnn
 \qed

 \begin{proposition}\label{ThmA07}
 	There exist constants $C>0$ and $\gamma\in(1,\alpha(1+\beta)/2)$  such that  for any $t\in[0,1]$ any $h>0$, and all $i\in\{1,2,4\}$ and $j\in\{5,7,8\}$,
 	\beqnn
	\mathbf{E}\big[|I_{T,i}(t,h)|^\alpha \cdot |I_{T,6}(t,h)|^\alpha\big] \leq Ch^{\gamma} \quad \mbox{and} \quad
	\mathbf{E}\big[|I_{T,3}(t,h)|^\alpha \cdot |I_{T,j}(t,h)|^\alpha\big]\leq Ch^{\gamma}.
 	\eeqnn
 \end{proposition}
 \proof 
 We just prove the result for $i=1$. 
 All other cases can be proved in the same way.
 As in the proof of Proposition~\ref{ThmA05}, there exists a constant $C>0$ such that
 \beqlb\label{eqnA.07}
 \mathbf{E}\big[|I_{T,1}(t,h)|^{2\alpha}\big] \leq C\frac{(Th)^\alpha+ Th}{T^\alpha}\quad \mbox{and}\quad \mathbf{E}\big[|I_{T,6}(t,h)|^{2\alpha}\big] \leq C|\delta_T|^\alpha.
 \eeqlb
 When $Th\geq 1$, by H\"older's inequality,
 \beqnn
 \mathbf{E}\big[|I_{T,1}(t,h)|^{\alpha}\cdot |I_{T,6}(t,h)|^{\alpha}\big]\leq C |\delta_T|^{\alpha/2} h^{\alpha/2}\leq CT^{-\beta\alpha/2}h^{\alpha/2}\leq Ch^{\alpha(1+\beta)/2}.
 \eeqnn
 Since $\beta>2/\alpha-1$, we see that (\ref{eqnA.04}) holds with $\gamma=\alpha(1+\beta)/2>1$. Moreover, $\gamma>1$ since $\alpha\beta>1$.
 When $Th<1$,  then
 \beqnn
 \mathbf{E}\big[|I_{T,1}(t,h)|^{2\alpha}\big]\leq CT^{1-\alpha}h.
 \eeqnn
 Applying Lemma~\ref{ThmA04}(2) to $I_{T,6}(t,h)$ with $t'=t+h-\delta_T$, $t''= t-h$, $f(t',t'',s,u)=   \int_{T(t-h-s)}^{T(t-s)} R(r,u)dr$, we have 
 \beqlb\label{eqnA.05.02}
 \lefteqn{\mathbf{E}[|I_{T,6}(t,h)|^{2\alpha}]
 \leq  \frac{C}{T^{\alpha-1}}\int_{t+h-\delta_T}^{t-h} \int_\mathbb{U} |\int_{T(t-h-s)}^{T(t-s)} R(r,u)dr|^{2\alpha}\nu_H(du)ds}\ar\ar \cr
 \ar\ar  + C  \int_{\mathbb{U}}   \nu_H(du) \int_{t+h-\delta_T}^{t-h} \Big( \int_{T(t-h-s)}^{T(t-s)} R(r,u)dr\Big)  ds 
 \cdot  \Big| \int_{t+h-\delta_T}^{t-h} \Big|\int_{T(t-h-s)}^{T(t-s)} R(r,u)dr\Big|^{\frac{2\alpha-1}{\alpha-1}}ds\Big|^{\alpha-1}.\quad
 \eeqlb
 From (\ref{eqn4.06}), we have 
 \beqnn
 \int_{T(t-h-s)}^{T(t-s)} R(r,u)dr\leq \frac{\|\phi(u)\|_\infty}{1-\|\phi_H\|_{L^1}}\cdot Th.
 \eeqnn
 From this and Condition~\ref{C1}, the first term on the right side of (\ref{eqnA.05.02}) can be bounded by 
 \beqnn
  \frac{C(Th)^{2\alpha}}{T^{\alpha-1}}\int_{t+h-\delta_T}^{t-h} \int_\mathbb{U} \|\phi(u)\|_\infty^{2\alpha}\nu_H(du)ds \leq CT^{\alpha+1}h^{2\alpha} \delta_T.
 \eeqnn
 Moreover, the second term also can be bounded by 
 \beqnn
 C (Th)^{2\alpha} |\delta_T|^\alpha \int_{\mathbb{U}}  \|\phi(u)\|_\infty^{2\alpha} \nu_0(du) \leq C (Th)^{2\alpha} |\delta_T|^\alpha.
 \eeqnn
 Putting all estimates above together, from the assumption that $Th<1$ we have
 \beqlb\label{eqnA.08}
 \mathbf{E}[|I_{T,6}(t,h)|^{2\alpha}] \ar\leq\ar  C (Th)^{2\alpha} |\delta_T|^\alpha.
 \eeqlb
 From H\"older's inequality and $\gamma\in(1,\alpha(1+\beta)/2)$, 
 \beqnn
 \mathbf{E}[|I_{T,1}(t,h)|^{\alpha}|I_{T,6}(t,h)|^{\alpha}]\leq C T^{1/2+\alpha(1-\beta)/2}  h^{1/2+\alpha} \leq C T^{1+\gamma-\alpha(1+\beta)/2}  h^{1+\gamma}\leq C h^{1+\gamma}.
 \eeqnn
 \qed

 \begin{corollary}\label{ThmA08}
 	There exist a constant $C>0$  such that for any $h<\delta_T/2$ and $T>0$,
 	\beqnn
 	\mathbf{E}[|I_{T,3}(t,h)|^{\alpha}|I_{T,6}(t,h)|^{\alpha}]\leq Ch^{\alpha\beta}.
 	\eeqnn
 \end{corollary}
 \proof From (\ref{eqnA.07}) and (\ref{eqnA.08}), we have
 \beqnn
 \mathbf{E}[|I_{T,3}(t,h)|^{2\alpha}]+\mathbf{E}[ |I_{T,6}(t,h)|^{2\alpha}] \leq
 \left\{  \begin{array}{ll}
 	C|\delta_T|^\alpha, & \mbox{if }Th\geq 1;\vspace{5pt}\\
 	C (Th)^{2\alpha} |\delta_T|^\alpha, & \mbox{if }Th< 1.
 \end{array}
 \right.
 \eeqnn
 By H\"older's inequality, when $Th\geq 1$,
 \beqnn
 \mathbf{E}[|I_{T,3}(t,h)|^\alpha|I_{T,6}(t,h)|^\alpha]\ar\leq\ar \Big( \mathbf{E}[|I_{T,3}(t,h)|^{2\alpha}]\mathbf{E}[ |I_{T,6}(t,h)|^{2\alpha}] \Big)^{1/2}\leq   C |\delta_T|^\alpha \leq C h^{\alpha\beta}
 \eeqnn
 and when $Th<1$,
 \beqnn
 \mathbf{E}[|I_{T,3}(t,h)|^\alpha|I_{T,6}(t,h)|^\alpha]\ar\leq\ar C T^{2\alpha-\beta \alpha} h^{2\alpha}\leq C (Th)^{2\alpha-\beta \alpha} h^{\alpha\beta}\leq Ch^{\alpha\beta}.
 \eeqnn
 \qed

 \subsection{Proof of Lemma~\ref{Thm405}}\label{AppendixB}

 In this section we give a detailed proof of the weak convergence of the sequence $\{\sqrt{T}\varepsilon^\psi_{T,H}(t):t\geq 0\}_{T\geq 1}$ to $0$ in the space $\mathbb{D}([0,\infty),\mathbb{R})$ under Condition~\ref{C1} and \ref{C2}.
 The weak convergence of $\{\sqrt{T}\varepsilon^\psi_{T,I}(t):t\geq 0\}_{T\geq 1}$ can be proved similarly. 
 
 Redefine $\delta_T:=T^{-\beta}$ with $\beta\in(\frac{1}{\alpha}, 1-\frac{1}{2\theta_1})$. We decompose the error process $\sqrt{T}\varepsilon^\psi_{T,H}(t)$ into the following two parts:
 \beqlb
 \varepsilon^\psi_{T,H,1}(t)\ar:=\ar  \frac{1}{\sqrt{T}}\int_0^{t-\delta_T} \int_{\mathbb{U}}  \psi^{\mathrm{c}}(T(t-s),u) N_H(dTs,du),\quad \label{eqnA.12}\\
 \varepsilon^\psi_{T,H,2}(t)\ar:=\ar  \frac{1}{\sqrt{T}}\int_{t-\delta_T}^{t} \int_{\mathbb{U}} \psi^{\mathrm{c}}(T(t-s),u) N_H(dTs,du).\label{eqnA.13}
 \eeqlb
 For the first term, we have 
 \beqlb
 |\varepsilon^\psi_{T,H,1}(t)|\ar\leq\ar \frac{1}{\sqrt{T}}\int_0^{t-\delta_T} \int_{\mathbb{U}}  \sup_{s\leq t-\delta_T }|\psi^{\mathrm{c}}(T(t-s),u)| N_H(dTs,du)\cr
 \ar\leq\ar  \frac{1}{\sqrt{T}}\int_0^{t-\delta_T} \int_{\mathbb{U}}  \sup_{s\geq T^{1-\beta}}|\psi^{\mathrm{c}}(s,u) | N_H(dTs,du).
 \eeqlb
 From Lemma~\ref{Thm201} and Condition~\ref{C2},
 \beqlb
 \mathbf{E}\Big[ \sup_{t\in[0,1]}  |\varepsilon^\psi_{T,H,1}(t)| \Big] \ar\leq\ar  \sqrt{T}\int_0^1 \mathbf{E}[ Z(Ts)]ds \int_{\mathbb{U}}  \sup_{s\geq T^{1-\beta}}|\psi^{\mathrm{c}}(s,u) | \nu_H(du)\leq  T^{\frac{1}{2} -\theta_1(1-\beta)} ,
 \eeqlb
 which vanishes as $T\to\infty$.  
 We now consider $\varepsilon^\psi_{T,H,2}(\cdot)$. For any $t\geq 0$, 
 \beqlb
 |\varepsilon^\psi_{T,H,2}(t)|\ar\leq \ar  \frac{1}{\sqrt{T}}\int_{t-\delta_T}^{t} \int_{\mathbb{U}} |\psi^{\mathrm{c}}(T(t-s),u) |N_H(dTs,du)\cr
 \ar\leq\ar  \frac{1}{\sqrt{T}}\int_{t-\delta_T}^{t} \int_{\mathbb{U}} \|\psi^{\mathrm{c}}(u) \|_\infty N_H(dTs,du).
 \eeqlb 
 Thus the proof would be finished if we can prove the following claim: for any $\epsilon>0$, we have as $T\to\infty$,
 \beqlb\label{eqnA.19}
 \mathbf{P}\Big\{\sup_{t\in[0,1]}  \frac{1}{\sqrt{T}}\int_{t-\delta_T}^{t} \int_{\mathbb{U}} \|\psi^{\mathrm{c}}(u) \|_\infty N_H(dTs,du)>\epsilon \Big\} \to 0.
 \eeqlb
 Indeed, from Chebyshev's inequality,  this probability can be bounded by
 \beqnn
 \lefteqn{ \mathbf{P}\Big\{ \sup_{i=0,\cdots [T^\beta]-1}\frac{1}{\sqrt{T}} \int_{iT^{-\beta}}^{(i+2)T^{-\beta}} \int_{\mathbb{U}}\|\psi^{\mathrm{c}}(u) \|_\infty N_H(dTs,du)>\epsilon \Big\}}\ar\ar\cr
 \ar\leq\ar  \sum_{i=0}^{ [T^\beta]-1}\mathbf{P}\Big\{ \frac{1}{\sqrt{T}} \int_{iT^{1-\beta}}^{(i+2)T^{1-\beta}} \int_{\mathbb{U}}\|\psi^{\mathrm{c}}(u) \|_\infty N_H(ds,du)>\epsilon \Big\}\cr
 \ar\leq\ar  \sum_{i=0}^{ [T^\beta]-1}\frac{1}{T^\alpha \epsilon^{2\alpha}} \mathbf{E}\Big[ \Big| \int_{iT^{1-\beta}}^{(i+2)T^{1-\beta}} \int_{\mathbb{U}}\|\psi^{\mathrm{c}}(u) \|_\infty N_H(ds,du)\Big|^{2\alpha} \Big].
 \eeqnn
 From the Cauchy-Schwarz inequality and (\ref{eqn2.04}), 
 \beqlb\label{eqnB.07}
 \lefteqn{\mathbf{E}\Big[ \Big| \int_{iT^{1-\beta}}^{(i+2)T^{1-\beta}} \int_{\mathbb{U}}\|\psi^{\mathrm{c}}(u) \|_\infty N_H(ds,du)\Big|^{2\alpha} \Big]}\ar\ar\cr
 \ar\leq\ar C\mathbf{E}\Big[ \Big| \int_{iT^{1-\beta}}^{(i+2)T^{1-\beta}} Z(s)ds\int_{\mathbb{U}}\|\psi^{\mathrm{c}}(u) \|_\infty \nu_H(du)\Big|^{2\alpha} \Big]\cr
 \ar\ar + C\mathbf{E}\Big[ \Big| \int_{iT^{-\beta}}^{(i+2)T^{-\beta}} \int_{\mathbb{U}} \int_0^{Z(Ts-)}\|\psi^{\mathrm{c}}(u) \|_\infty \tilde{N}_0(dTs,du,dz)\Big|^{2\alpha} \Big].
 \eeqlb
 From Condition~\ref{C2}, H\"older's ineqaulity and Lemma~\ref{ThmA04}(1), we see that the first term on the right side of this inequality can be bounded by
 \beqlb\label{eqnB.08}
 C\mathbf{E}\Big[ \Big| \int_{iT^{1-\beta}}^{(i+2)T^{1-\beta}} Z(s)ds\Big|^{2\alpha} \Big]
 \leq C T^{(2\alpha-1)(1-\beta)} \int_{iT^{1-\beta}}^{(i+2)T^{1-\beta}}\mathbf{E}[ |Z(s)|^{2\alpha} ] ds\leq CT^{2\alpha(1-\beta)}.
 \eeqlb
 Applying Lemma~\ref{ThmA04}(2) to the second term on the right side of  (\ref{eqnB.07}) with $t'= i\cdot T^{-\beta}$, $t''=(i+2)T^{-\beta}$ and $f(t',t'',s,u)= \sup_{r\geq 0} |\psi^{\mathrm{c}}(r,u) |$,  from Condition~\ref{C2} we have 
 \beqnn
 \lefteqn{\mathbf{E}\Big[ \Big| \int_{i\cdot T^{-\beta}}^{(i+2)T^{-\beta}} \int_{\mathbb{U}} \int_0^{Z(Ts-)}\|\psi^{\mathrm{c}}(u) \|_\infty \tilde{N}_0(dTs,du,dz)\Big|^{2\alpha} \Big]}\ar\ar\cr
 \ar\leq\ar C(T^{1-\beta}  +T^{\alpha(1-\beta) })\int_{\mathbb{U}}   \|\psi^{\mathrm{c}}(u) \|_\infty^{2\alpha} \nu_H(du) \leq CT^{\alpha(1-\beta) }.
 \eeqnn
 Putting all estimates above together, we have 
 \beqnn
 \mathbf{P}\Big\{\sup_{t\in[0,1]}  |\varepsilon^\psi_{T,H,2}(t)|>\epsilon \Big\} \leq  \frac{ C}{ \epsilon^{2\alpha}}T^{\alpha   -(2\alpha-1)\beta}.
 \eeqnn
 which vanishes as $T\to\infty$ since $\beta>\frac{\alpha}{2\alpha-1}$. Here we have gotten (\ref{eqnA.19}) and the proof of Lemma~\ref{Thm405} has been finished.

 %
 %

 \section{Application to budding microbes in a host}
 \setcounter{equation}{0}
 \medskip
 
 In this section, we apply our limit theorems for marked Hawkes point measures to study the asymptotics of the amounts of toxins released by budding microbes in a host. 
 Let $X(0)\in\mathbb{N}$ be the number of microbes in the host at time $0$. 
 Different to binary fission where the fully grown parent cell splits into two equally sized daughter cells, small buds usually form at one end of mother cell or on filaments called prosthecae at random budding times $0<\tau_1<\tau_2<\cdots $ and separate as new microbes when mature.
 The mother cell produces buds before dying, spreading out of the host or being killed by the host.
 We assume the \textit{life-length} of microbes from birth to death is randomly distributed with common probability law $\Lambda_H(dt)$ defined on $(0,\infty)$.
 The life-length $\Lambda_H(dt)$ is rarely exponentially distributed; see \cite[Table~4]{HolbrookMenninger2002} and \cite[Figure 2-4]{Wood2004}.
 
 Conditioned on the life-length $y$ and age $t$, we assume that the mother microbe produces new buds at the \textit{budding rate} $\gamma_H(t,y)$, where $\gamma_H(
 \cdot)$ is a nonnegative funciton on $\mathbb{R}_+^2$.
 Usually, the budding rate is low during the growth stage.
 After separating from the mother cell, the budding rate increases to its highest level.
 However, as bud scars accumulate on its surface, the microbe enters into the senescence stage and the budding rate starts to decrease; see \cite[Figure~2]{Jiang2000}. Without loss of generality, we assume that the budding rate $\gamma_H(\cdot)$ is bounded.
 
 Usually, only one bud forms on the mother cell at each budding time. However, multiple-budding also happen in the reproduction of budding viruses such as HIV; see \cite[p.416]{Tortora2016}.
 Hence we assume that a random number of buds forms on the mother cell at each budding time according to the probability law $\mathbf{p}_H=(p_{H,1},p_{H,2},\cdots)$ with  generating function $\{ g_H(z):=\sum_{k=1}^\infty p_{H,k} z^k : z\in[0,1] \}$.
 When $p_{H,1}=1$ and $\Lambda(dt)$ is an exponential distribution, the budding reproduction reduces to binary fission\footnote{Ackermann et al. \cite{Ackermann2003} and Stephens \cite{Stephens2005} presented evidence for the existence of asymmetric division in binary fission, where the two daughter microbes are not equal.
 Asymmetric binary fission is captured by our model.}. 
 In addition to budding, microbes may immigrate from external sources or neighbouring hosts at random times $0<\sigma_1<\sigma_2<\cdots$. 
 To simplify the analysis, we assume that the arrivals of immigrating microbes follow some Poisson point process with intensity $\lambda_1=1$ and that the number of invading bacteria is distributed according to the probability distribution $\mathbf{p}_I=(p_{I,1},p_{1,2},\cdots)$ with generating function $\{ g_I(z):=\sum_{k=1}^\infty p_{I,k} z^k : z\in[0,1] \}$.
%
We consider the ancestors at time $0$ as the $0$-th immigration and allow the budding rate function $\gamma_I(\cdot)$ of immigrating microbes to be different from $\gamma_H(\cdot)$. 
 
 According to their  origins, we classify the microbes in the host into the following three classes:
 \begin{enumerate}
 	\item[] $\mathcal{I}_0$: Ancestors at time $0$;
 	
 	\item[] $\mathcal{I}_i$:  Bacteria migrating into the host at the $i$-th immigrating time;
 	
 	\item[] $\mathcal{B}_i$:  Buds produced at the $i$-th budding time.
 \end{enumerate}
  Let $\mathbf{B}(t)$ be the total budding rate of all bacteria alive at time $t$. It can be written as
 \beqlb\label{eqn6.14}
 \mathbf{B}(t)\ar=\ar \sum_{j\in\mathcal{I}_0} \gamma_I(t, \ell_{0,j})+\sum_{\sigma_i\leq t} \sum_{j\in\mathcal{I}_i} \gamma_I(t-\sigma_i, \ell_{i,j})  +\sum_{\tau_i\leq t} \sum_{j\in\mathcal{B}_i} \gamma_H(t-\tau_j, \ell_{i,j}),
 \eeqlb
 where $\ell_{i,j}$ is the life-length of $j$-th microbes in $\mathcal{B}_i$ or $\mathcal{I}_i$.
 Here the first sum on the right side of the above equality is the total budding rate of the ancestors in the host at time $t$,
 the inner-sum in the second term is the total budding rate of the bacteria immigrating into the host during the $j$-th invasion, and
 the inner-sum in the third term is the total budding rate of all buds formed at the $j$-th budding time. Microbes do not only  produce offsprings but also infect the host by releasing toxins and attacking the host cell.
 For instance,  \textit{Candida albicans} (C. albicans) in the gastrointestinal and genitourinary tract do not only release a toxin called \textit{Candidalysin} but also alkalinize phagosomal by physical rupture.
 Denote by $\mathtt{T}(t,y)$ the \textit{cumulative toxins} released or \textit{cumulative damage} made by a microbe with life-length $y$ to the host up to age $t$.
 After dying or being killed by the host, the microbe stops releasing toxins, i.e. $\mathtt{T}(t,y)=\mathtt{T}(y,y)$ when $t\geq y$ .
 Because of the diversity in bacteria, their toxin release functions are usually different.
 Let $\mathbb{T}$ be the collection of toxin cumulative functions:
 \beqlb\label{eqn6.03}
 \mathbb{T}:=\big\{ \mathtt{T}\ar: \ar \mathtt{T}(t,y) \mbox{ is a nondecreasing function with } \mathtt{T}(t,y)= \mathtt{T}(y,y) \mbox{ if } t\geq y \big\}.
 \eeqlb
 For any $\mathtt{T}\in\mathbb{T}$, denote by $\mathtt{T}^{\mathrm{c}}(t,y)$ the $\textit{unreleased toxins}$ of a microbe with life-length $y$ at age $t$, i.e. $\mathtt{T}^{\mathrm{c}}(t):=\mathtt{T}(y,y)-\mathtt{T}(t,y)$.
 We assume that the toxin function of microbes born in the hosts is distributed according to the law $m_H (d\mathtt{T})$  and  the toxin function of microbes immigrating into the host is distributed according to the law $m_I (d\mathtt{T})$. We also assume that each microbe picks up its toxin function independently. 
 
 Most microbes release toxins continuously during their life. For instance, C.albicans release toxins continuously during hyphal formation.
 In this case, conditioned on the life-length $y$ and age $t$, we may assume that the bacteria releases toxins at rate $\varphi(t)$ and $\mathtt{T}(t,y) := \int_0^{t\wedge y} \varphi(s)ds$, where $\{\varphi(s) :s\geq 0 \}$ is a functional-valued random variable. 
 Other microbes release toxins immediately when they decompose. In this case  $\mathtt{T}(t,y) = \vartheta\cdot \mathbf{1}_{\{t\geq y  \}}$, where $\vartheta$, a $\mathbb{R}_+$-valued random variable, is the amount of toxin released by the microbe at the time of death.
 For any $t\geq 0$, let $\mathbf{T}(t)$ be the \textit{total cumulative toxins} released by the entire population up to time $t$. Similar to the representation of $\mathbf{B}(t)$, we can  represent $\mathbf{T}(t)$  as  
 \beqlb\label{eqn6.04}
 \mathbf{T}(t):=  \sum_{j\in\mathcal{I}_0} \mathtt{T}_{j}(t, \ell_{0,j})+\sum_{\sigma_i\leq t} \sum_{j\in\mathcal{I}_i}\mathtt{T}_{j}(t-\sigma_i, \ell_{i,j})  +\sum_{\tau_i\leq t} \sum_{j\in\mathcal{B}_i}\mathtt{T}_{j}(t-\tau_j, \ell_{i,j}).
 \eeqlb
 

 We are now going to describe the budding and population dynamics in terms of a marked Hawkes process. To this end, we choose the mark space $\mathbb{U} := \{H,I \}\times \mathbb{Z}_+\times\mathbb{R}_+^{\mathbb{Z}_+}\times \mathbb{T}^{\mathbb{Z}_+}$. 
 Here $u=(i,k,\boldsymbol{y},\boldsymbol{T})\in \mathbb{U}$ means that $i\in\{ H,I\}$, $k\in\mathbb{Z}_+$,  $\boldsymbol{y}:=(y_1,\cdots,y_k)\in\mathbb{R}_+^k$ and $\boldsymbol{T}:=(\mathtt{T}_1,\cdots,\mathtt{T}_k)\in  \mathbb{T}^k$.
 We record the information of buddings and invasions  with two sequences of i.i.d.  $\mathbb{U}$-valued random variables  $\{\xi_j: j=1,2,\cdots\}$ and $\{\eta_j: j=1,2,\cdots\}$ respectively, i.e., $\xi_j/\eta_j=(H/I, k,\boldsymbol{y},\boldsymbol{T})$ means that there are $k$ buds/microbes with life-length $(y_1,\cdots,y_k )$ and toxin cumulative function $ (\mathtt{T}_1,\cdots,\mathtt{T}_k)$ splitting from the mother cell/immigrating at time $\tau_j/\sigma_j$.
  According to our previous assumption, $\xi$ and $\eta$ have the probability laws $\nu_H(du)$ and $\nu_I(du)$ defined as follows: for $i'\in\{H,I \}$,
  \beqlb\label{eqn6.17}
  \nu_{i'}(du):=\nu_{i'}(di,dk,d\boldsymbol{y},d\boldsymbol{T})\ar=\ar \mathbf{1}_{i'}(di) \sum_{j=1}^\infty p_{i',j}\mathbf{1}_{\{j\}}(dk)\prod_{l=1}^j \Lambda_{i'}(dy_l)m_{i'}(d\mathtt{T}_l).
  \eeqlb
  As argued in Section~2, the total budding rate $\mathbf{B}(t)$ and the  total cumulative toxins  $\mathbf{T}(t)$ at time  $t$ can be represented as:
  \beqnn
  \mathbf{B}(t)\ar=\ar \sum_{i=1}^{X(0)} \gamma_I(t, \ell_{0,i}) +  \int_0^t \int_{\mathbb{U}}  \sum_{j=1}^k\gamma_i(t-s,y_j)N_I(ds,du) \cr
  \ar\ar + \int_0^t \int_{\mathbb{U}}\int_0^{\mathbf{B}(s-)}  \sum_{j=1}^k\gamma_i(t-s,y_j)N_0(ds,du,dz),\cr 
  \mathbf{T}(t)\ar=\ar \sum_{i=1}^{X(0)} \mathtt{T}_i(t, \ell_{0,i}) + \int_0^t \int_{\mathbb{U}}   \sum_{j=1}^k\mathtt{T}_j(t-s,y_j)N_I(ds,du)\cr
  \ar\ar  + \int_0^t \int_{\mathbb{U}}\int_0^{\mathbf{B}(s-)} \sum_{j=1}^k\mathtt{T}_j(t-s,y_j)N_0(ds,du,dz).
  \eeqnn
  where $N_0(ds,du,dz)$ and $N_I(ds,du)$ are two independent time-homogeneous Poisson random measures defined as before. 
 In this case, Condition~\ref{C1} and \ref{C2} reduce to the following condition.

   \begin{condition} \label{C3}
  	Recall the constants $\alpha$, $\theta_0$ and $\theta_1$. We assume that for $i\in\{H,I\}$, 
  	\beqlb\label{eqn6.05}
  	\int_0^\infty\Big[  \int_0^\infty t^{\theta_0}\gamma_i(t,y)dt+y^{2\alpha}\Big]  \Lambda_i(dy)+  \sum_{k=1}^\infty k^{2\alpha} p_{i,k} <\infty
  	\eeqlb
  	and
  	\beqlb\label{eqn6.06}
  	\sup_{t\geq 0} 	\int_0^\infty \Lambda_i(dy) \int_{\mathbb{T}}\Big[|\mathtt{T}(t,y)|^{2\alpha} + t^{\theta_1}   \mathtt{T}^{\mathrm{c}}(t,y)\Big] m_i(d\mathtt{T}) <\infty.
  	\eeqlb
  	
  \end{condition}
  Under this condition, the following quantities are well defined. For $i\in\{H,I \}$ ,  let $ g'_i(1):=\sum_{k=1}^\infty k p_{i,k}$  and $g''_i(1):= \sum_{k=1}^\infty k(k-1) p_{i,k}$.
  Moreover, for any $\kappa=1,2$, define
  \beqlb\label{eqn6.08}
  \|\gamma\|_{\Lambda^\kappa_i}^\kappa :=  \int_0^\infty \|\gamma_i(y)\|_{L^1}^\kappa \Lambda_i(dy),
  \quad \|\mathtt{T} \|_{\Lambda^\kappa_i}^\kappa:= \int_0^\infty\Lambda_i(dy)  \int_{\mathbb{T}} |\mathtt{T}(y,y)|^\kappa m_i(d\mathtt{T})
  \eeqlb
  and
  \beqlb\label{eqn6.09}
  \langle  \gamma, \mathtt{T}  \rangle_i :=  \int_0^\infty  \|\gamma_i(y)\|_{L^1} \Lambda_i(dy) \int_{\mathbb{T}} \mathtt{T}(y,y)m_i(d\mathtt{T}).
  \eeqlb
  From Lemma~\ref{Thm201}, we can see that $\sup_{t\geq 0}\mathbf{E}[\mathbf{B}(t)]<\infty$ if and only if $g'_H(1)\|\gamma\|_{\Lambda_H^1}<1$.
  Applying Lemma~\ref{Thm401}, Theorem~\ref{Thm404} and \ref{Thm406} with 
 \beqnn
 \phi(t,u)=  \sum_{j=1}^k\gamma_{i}(t,y_j) \quad \mbox{and}\quad \psi(t,u)=   \sum_{j=1}^k\mathtt{T}_j(t,y_j),
 \eeqnn
 we can get the following functional central limit theorem for $\{(\int_0^t \mathbf{B}(s)ds,\mathbf{T}(t)):t\geq 0  \}$.

 \begin{theorem}\label{Thm601}
 Suppose that $g'_H(1)\|\gamma\|_{L^1_\Lambda}<1$ and Condition~\ref{C3} holds. Then
 	\beqlb\label{eqn6.10}
 	\lefteqn{	\sqrt{T}\left( \begin{array}{c}
 	\frac{\int_0^{Tt}\mathbf{B}(s)ds}{T} -\frac{g'_I(1)\|\gamma\|_{\Lambda_I^1}}{1-g'_H(1)\|\gamma\|_{\Lambda_H^1}}\cdot t \vspace{5pt}\cr
 		\frac{\mathbf{T}(Tt)}{T}-g'_I(1) \frac{ \|\mathtt{T}\|_{\Lambda_I^1}+ g'_H(1)[\|\gamma\|_{\Lambda_I^1}\|\mathtt{T}\|_{\Lambda_H^1}-\|\gamma\|_{\Lambda_H^1 }\|\mathtt{T}\|_{\Lambda_I^1} ]}{1- g'_H(1) \|\gamma\|_{\Lambda_H^1}}\cdot t
 	\end{array}\right)}\ar\ar\cr
 \ar\to\ar  \sum_{i\in\{H,I \}} \int_0^t \int_{\mathbb{Z}_+}\int_{\mathbb{R}_+^{\mathbb{Z}_+}}\int_{\mathbb{T}^{\mathbb{Z}_+}}  \Phi_i(s,k,\boldsymbol{y},\boldsymbol{T}) W_i(ds,dk,d\boldsymbol{y},d\boldsymbol{T})
 	\eeqlb
 	weakly in the space $\mathbb{D}([0,\infty),\mathbb{R})$, where
 		\beqnn
 	\Phi^{\rm T}_i(s,k,\boldsymbol{y},\boldsymbol{T})= \Big(  \frac{\sum_{j=1}^k \|\gamma_i(y_j)\|_{L^1} }{1- g'_H(1)\|\gamma\|_{\Lambda_H^1}} \ , \  \sum_{j=1}^k \mathtt{T}_j(y_j,y_j)+ g'_H(1)\|\mathtt{T}\|_{\Lambda_H^1}  \frac{\sum_{j=1}^k \|\gamma_i(y_j)\|_{L^1} }{1- g'_H(1)\|\gamma\|_{\Lambda_H^1}} \Big),
 	\eeqnn
   and $W_H(ds,dk,d\boldsymbol{y},d\boldsymbol{T})$ and $W_I(ds,dk,d\boldsymbol{y},d\boldsymbol{T})$ are two independent Gaussian white noises on $(0,\infty)\times  \mathbb{Z}_+ \times \mathbb{R}_+^{\mathbb{Z}_+}\times \mathbb{T}^{\mathbb{Z}_+}$ with intensities
 	\beqnn
 	\frac{g'_I(1)\|\gamma\|_{\Lambda_I^1}}{1- g'_H(1)\|\gamma\|_{\Lambda_H^1}}ds\nu_H(\{H \}, dk,d\boldsymbol{y},d\boldsymbol{T}) \quad \mbox{and}\quad
 	ds\nu_I(\{I\}, dk,d\boldsymbol{y},d\boldsymbol{T}).
 	\eeqnn

 	 \end{theorem}
 	
 	From the property of stochastic integral with respect to Gaussian white noise, we can see that the limit process in (\ref{eqn6.10}) is a linear combination of two independent two-dimensional Brownian motions, i.e.
 	\beqnn
 	\int_0^t \int_{\mathbb{Z}_+}\int_{\mathbb{R}_+^{\mathbb{Z}_+}}\int_{\mathbb{T}^{\mathbb{Z}_+}}  \Phi_H(s,k,\boldsymbol{y},\boldsymbol{T}) W_H(ds,dk,d\boldsymbol{y},d\boldsymbol{T})= \Big|\frac{g'_I(1)\|\gamma\|_{\Lambda_I^1}}{1- g'_H(1)\|\gamma\|_{\Lambda_H^1}} \Big|^{1/2}\cdot \left( \begin{array}{l}  c_{H1} B_{H1}(t) \\  c_{H2} B_{H2}(t)  \end{array}\right)
 	\eeqnn
 	and
 		\beqnn
 	\int_0^t \int_{\mathbb{Z}_+}\int_{\mathbb{R}_+^{\mathbb{Z}_+}}\int_{\mathbb{T}^{\mathbb{Z}_+}}  \Phi_I(s,k,\boldsymbol{y},\boldsymbol{T}) W_I(ds,dk,d\boldsymbol{y},d\boldsymbol{T})=  \left( \begin{array}{l}  c_{I1} B_{I1}(t) \\  c_{I2} B_{I2}(t)  \end{array}\right),
 	\eeqnn
 	 where for $i\in\{ H,I \}$,
 	 \beqnn
 	 |c_{i1}|^2\ar=\ar \frac{g'_i(1) \|\gamma\|^2_{\Lambda_i^2} + g''_i(1) \|\gamma\|^2_{\Lambda_i^1}}{|1- g'_H(1)\|\gamma\|_{\Lambda_H^1}|^2} , \cr
 	 |c_{i2}|^2\ar=\ar g'_i(1)\Big|\|\mathtt{T}\|_{\Lambda_i^2}+\frac{g'_H(1)\|\mathtt{T}\|_{\Lambda_H^1} \|\gamma\|_{\Lambda_i^2} }{1- g'_H(1)\|\gamma\|_{\Lambda_H^1}}\Big|^2  +g''_i(1)\Big|\|\mathtt{T}\|_{\Lambda_i^1}+\frac{g'_H(1)\|\mathtt{T}\|_{\Lambda_H^1} \|\gamma\|_{\Lambda_i^1} }{1- g'_H(1)\|\gamma\|_{\Lambda_H^1}}\Big|^2 \cr
 	 \ar\ar
 	 +  \frac{2g'_H(1)\|\mathtt{T}\|_{\Lambda_H^1}g'_i(1) }{1- g'_H(1)\|\gamma\|_{\Lambda_H^1}} \Big(\langle \gamma,\mathtt{T}  \rangle_i 
 	 - \|\mathtt{T}\|_{\Lambda_i^2}\|\gamma\|_{\Lambda_i^2}\Big)
 	 \eeqnn
 	 and $(B_{H1}(t),B_{H2}(t))$ and $(B_{I1}(t),B_{I2}(t))$ are two independent two-dimensional Brownian motions  with 
 	 \beqnn 
 	 \langle  B_{i1},B_{i2}  \rangle_t \ar=\ar \frac{t}{c_{i1}c_{i2}}\Big( \frac{g'_i(1)\langle \gamma,\mathtt{T} \rangle_i + g''_i(1) \|\gamma\|_{\Lambda_i^1}\|\mathtt{T}\|_{\Lambda_i^1}}{ 1- g'_H(1)\|\gamma\|_{\Lambda_H^1}}+\frac{g'_H(1)\|\mathtt{T}\|_{\Lambda_H^1} ( g'_i(1) \|\gamma\|_{\Lambda_i^2}^2 + g''_i(1)\|\gamma\|_{\Lambda_i^1}^2  ) }{|1- g'_H(1)\|\gamma\|_{\Lambda_H^1}|^2} \Big).
 	 \eeqnn

 We close this section with two special and important cases. 
 The cumulative demage made by budding viruses (e.g. HIV) to the host is usually measured by the total number of viruses.
 In this case, we may assume that the toxin cumulative process satisfies that $\mathtt{T}(t,y)\equiv 1$ with probability one.  
 In this case, Condition~\ref{C3} holds in this case and the total toxin cumulative is given by $\mathbf{T}(t)= |\mathcal{I}(t)|$, the \textit{total progeny} of the microbial  population.
 
 \begin{proposition}\label{Thm602}
 	Suppose that $g'_H(1)\|\gamma\|_{L^1_\Lambda}<1$ and Condition~\ref{C3} holds, we have as $T\to\infty$,
 	\beqlb\label{eqn6.18}
  \lefteqn{	\sqrt{T} \Big(\frac{1}{T} |\mathcal{I}(Tt)|-g'_I(1) \frac{ 1+ g'_H(1)[\|\gamma\|_{\Lambda_I^1} -\|\gamma\|_{\Lambda_H^1}   ]}{1- g'_H(1) \|\gamma\|_{\Lambda_H^1}}\cdot t\Big)}\ar\ar\cr
  \ar\ar \to \Big|\frac{g'_I(1)\|\gamma\|_{\Lambda_I^1}}{1- g'_H(1)\|\gamma\|_{\Lambda_H^1}} \Big|^{1/2} \cdot c_{H2} B_{H2} (t) + c_{I2} B_{I2}(t),
 	\eeqlb
 	weakly in the space $\mathbb{D}([0,\infty),\mathbb{R})$, where  
 	 	\beqlb\label{eqn6.19}
 	|c_{i2}|^2\ar=\ar g'_i(1)\Big|1+\frac{g'_H(1) \|\gamma\|_{\Lambda_i^2} }{1- g'_H(1)\|\gamma\|_{\Lambda_H^1}}\Big|^2  +g''_i(1)\Big|1+\frac{g'_H(1) \|\gamma\|_{\Lambda_i^1} }{1- g'_H(1)\|\gamma\|_{\Lambda_H^1}}\Big|^2 \cr
 	\ar\ar
 	+  \frac{2g'_H(1) g'_i(1) }{1- g'_H(1)\|\gamma\|_{\Lambda_H^1}} \Big( \|\gamma\|_{\Lambda_i^1}
 	-  \|\gamma\|_{\Lambda_i^2}\Big).
 	\eeqlb
 \end{proposition}
 
 As we have mentioned before, most microbes release toxins continuously during their life. 
 Specially, assume that $\varphi(t)\equiv 1$ and $\mathtt{T}(t,y)=t\wedge y$,
 the total toxin cumulative process is 
 \beqlb\label{eqn6.20}
 \mathbf{T}(t)\ar=\ar   \int_0^t X(s)ds,
 \eeqlb
 where $X(t)$ denote the total number of microbes alive at time $t$ in the host.
 In this case, $\mathbf{T}(\cdot)$ is  usually called the \textit{integral of population}. 
In this case, Condition~\ref{C3} holds when $\alpha>3/2$ and the functional central limit theorem holds for the integral of population. 
The central limit theorem for the integral of population was proved in \cite{Pakes1975}.
 
 \begin{proposition}\label{Thm603}
 	 Suppose that $g'_H(1)\|\gamma\|_{L^1_\Lambda}<1$ and Condition~\ref{C3} holds with $\alpha>3/2$. Let $ \|\Lambda_i^\kappa\|^\kappa:= \int_0^\infty y^k\Lambda_i(dy)$ for $\kappa=1,2$ and $i\in\{H,I \}$. we have as $T\to\infty$,
 	\beqlb\label{eqn6.22}
 	\ar\ar \sqrt{T} \Big( \frac{1}{T} \int_0^{Tt} X(s)ds-g'_I(1) \frac{  \|\Lambda_I^1\| + g'_H(1)[\|\gamma\|_{\Lambda_I^1} \|\Lambda_H^1\|-\|\gamma\|_{\Lambda_H^1 } \|\Lambda_I^1 \|]}{1- g'_H(1) \|\gamma\|_{\Lambda_H^1}}\cdot t\Big) \cr
 	\ar\ar  \to  \Big|\frac{g'_I(1)\|\gamma\|_{\Lambda_I^1}}{1- g'_H(1)\|\gamma\|_{\Lambda_H^1}} \Big|^{1/2}  \cdot c_{H2} B_{H2}  (t)+ c_{I2} B_{I2} (t),
 	\eeqlb
 	weakly in the space $\mathbb{D}([0,\infty),\mathbb{R})$, where 
 	\beqlb\label{eqn6.23}
 	|c_{i2}|^2\ar=\ar g'_i(1)\Big|\|\Lambda_i^2\|+\frac{g'_H(1)\|\Lambda_H^1\| \|\gamma\|_{\Lambda_i^2} }{1- g'_H(1)\|\gamma\|_{\Lambda_H^1}}\Big|^2  +g''_i(1)\Big|\|\Lambda_i^1\|+\frac{g'_H(1)\|\Lambda_H^1\| \|\gamma\|_{\Lambda_i^1} }{1- g'_H(1)\|\gamma\|_{\Lambda_H^1}}\Big|^2 \cr
 	\ar\ar
 	+  \frac{2g'_H(1)\|\Lambda_H^1\|g'_i(1) }{1- g'_H(1)\|\gamma\|_{\Lambda_H^1}} \Big(\int_0^\infty y\|\gamma_i(y)\|_{L^1}\Lambda_i(dy)  
 	-\|\Lambda_i^2\|\|\gamma\|_{\Lambda_i^2}\Big).
 	\eeqlb
 	
 \end{proposition}
 
 \section{Conclusion}
 \setcounter{equation}{0}
 \medskip
 
 This paper established functional limit theorems for marked Hawkes point measures with immigration and their shot noise processes. 
 We proved that a suitably normalized point measure converges in distribution to the sum of a Gaussian white noise and a lifting maps of a Brownian motion. 
 The Brownian motion results from the cumulative intensity process that can be viewed as a form of common factor for the arrival of events with different marks. 
 Our limit theorems were used to analyze the population dynamics of microbes in a host and its interaction with that host. 
 At least three interesting problems were left open. 
 
 First, we assumed that the impact of each event to the arrivals of future events is short-term; see the first integral in (\ref{eqn3.07}). 
 For the long memory case, we expect a limiting diffusion in terms of a Gaussian white noise and a lifting of a continuous Gaussian process, which can be decompose into a Brownian motion and its Holmgren-Riemann-Liouville integral when the kenel $\phi(t,u)$ is regularly varying in time. Second, we considered the case that the distribution of marks is light-tailed; see the second integral in (\ref{eqn3.07}). If this integral is infinite for $\alpha=1$, we expect the weak convergence of the time-spatial rescaled marked Hawkes point measure to a lifting map of some nonnegative jump-diffusion process. Third, it would be interesting to derive a large deviation principle for marked Hawkes point measures in any cases above. These three problems will be addressed in future research.

\bibliographystyle{siam}

\bibliography{ReferenceforMHP}

\end{document}